\documentclass[11pt,reqno,oneside]{amsart}
\usepackage{verbatim,amsmath,amssymb,cite,epsfig,xspace,color,subfigure}
\usepackage[foot]{amsaddr}

\numberwithin{equation}{section}

\theoremstyle{plain}
\newtheorem{theorem}{Theorem}[section]

\newtheorem{corollary}{Corollary}[section]
\theoremstyle{definition}

\theoremstyle{remark}
\newtheorem{remark}{Remark}[section]

\newcommand{\ignore}[1]{}

\newcommand{\la}{\langle}
\newcommand{\ra}{\rangle}

\newcommand{\eps}{\epsilon}
\renewcommand{\c}{{\rm c}}
\newcommand{\const}{{\rm const}}

\definecolor{AlexColor}{rgb}{0.0,0,0}
\definecolor{HelenColor}{rgb}{0.0,0,0}
\definecolor{MitchColor}{rgb}{0.0,0,0}
\definecolor{ChristophColor}{rgb}{0.0,0,0}
\newcommand{\helen}[1]{{\color{HelenColor}#1}}
\newcommand{\as}[1]{{\color{AlexColor}#1}}

\usepackage{ulem}

\usepackage{cancel}

\newcommand{\co}[1]{{\color{ChristophColor}#1}}

\begin{document}

\title
[Theory-based Benchmarking of the B-QCF method.]
{Theory-based Benchmarking of the Blended Force-Based Quasicontinuum Method}

\author{Xingjie Li}
\address{Xingjie Li\\ 182 George St. \\
  Providence \\ RI 02912 \\ USA\\xingjie li@brown.edu}

\author{Mitchell Luskin}
\address{M. Luskin (Corresponding Author)\\ 127 Vincent Hall \\ 206 Church St. SE \\
  Minneapolis \\ MN 55455 \\ USA\\ luskin@umn.edu\\ Phone: 612-625-6565\\ FAX 612-626-2017}

\author{Christoph Ortner}
\address{C. Ortner\\ Mathematics Institute \\ Zeeman Building \\
  University of Warwick \\ Coventry CV4 7AL \\ UK\\christoph.ortner@warwick.ac.uk}

\author{Alexander V. Shapeev}
\address{A. V. Shapeev\\ 127 Vincent Hall \\ 206 Church St. SE \\
  Minneapolis \\ MN 55455 \\ USA\\alexander@shapeev.com}


\thanks{ This work was supported in part by the NSF PIRE Grant
  OISE-0967140, DOE Award DE-SC0002085, and AFOSR Award
  FA9550-12-1-0187. CO was supported by EPSRC grant EP/H003096
  ``Analysis of atomistic-to-continuum coupling methods.''}



\date{\today}

\begin{abstract}
We formulate an atomistic-to-continuum coupling method based on
blending atomistic and continuum forces. Our precise choice of
blending mechanism is informed by theoretical predictions. We
present a range of numerical experiments studying the accuracy
of the scheme, focusing in particular on its stability. These
experiments confirm and extend the theoretical predictions, and
demonstrate a superior accuracy of B-QCF over energy-based
blending schemes.
\end{abstract}

\maketitle{\thispagestyle{empty} 
\section{Introduction}
Atomistic-to-continuum coupling methods (a/c methods) have been
proposed to increase the computational efficiency of atomistic
computations involving the interaction between local crystal defects
with long-range elastic
fields~\cite{curt03,LinP:2006a,Miller:2003a,Shimokawa:2004,E:2004,Miller:2008,Legoll:2005,bqce11};
see \cite{atc.acta} for a recent review of a/c coupling methods and their numerical analysis.
Energy-based methods in this class, such as the quasicontinuum model
(denoted QCE \cite{Ortiz:1995a}), exhibit spurious interfacial forces
(``ghost forces'') even under uniform strain~\cite{Shenoy:1999a,
  Dobson:2008a}.  The effect of the ghost force on the error in
computing the deformation and the lattice stability by the QCE
approximation has been analyzed in~\cite{Dobson:2008a, Dobson:2008c,
  mingyang, doblusort:qce.stab}, {\helen{where lattice stability refers
  to the positive
  definiteness of the Hessian matrix of the total potential energy}}. The development of more accurate
energy-based a/c methods is an ongoing process
\cite{Shimokawa:2004,E:2004, Shapeev2D:2011,Shapeev3D:2011,
  LuskinXingjie.qnl1d,OrtnerZhang:2011,XiaoBely:2002}.

An alternative approach to a/c coupling is the force-based
quasicontinuum (QCF) approximation~\cite{doblusort:qcf.stab, qcf.stab,
  curt03, Miller:2003a, Lu.bqcf:2011}, but the non-conservative and
indefinite equilibrium equations make the iterative solution and the
determination of lattice stability more
challenging~\cite{qcf.iterative, qcf.stab, DobShapOrt:2011}. Indeed,
it is an open problem whether the (sharp-interface) QCF method is
stable in dimension greater than one. Although some recent results
  in this direction exist \cite{LuMing12}, it is still
  unclear to what extent they can be extended for general atomistic
  domains and in the presence of defects.

Many blended a/c coupling methods have been proposed in the
literature, e.g.,
~\cite{xiao:bridgingdomain,badia:onAtCcouplingbyblending,bridging,badia:forcebasedAtCcoupling,seleson:bridgingmethods,fish:concurrentAtCcoupling,prudhomme:modelingerrorArlequin,bauman:applicationofArlequin,XiBe:2004}.
In \cite{BQCF}, we formulated a blended force-based quasicontinuum
(B-QCF) method, similar to the method proposed in \cite{Lu.bqcf:2011},
which smoothly blends the forces of the atomistic and continuum model
instead of the sharp transition in the QCF method.  Under the
  simplifying assumption that deformation is homogeneous, we
  established sharp conditions under which a linearized B-QCF operator
  is positive definite, which effectively guarantees stability of the
  numerical scheme.
  \as{Surprisingly, the required blending width to ensure positive
  definiteness of the linearized B-QCF operator is \emph{asymptotically} small (however typical prefactors in the relative size of the blending region are not predicted by the theory).
  The one-dimensional theory developed in \cite{BQCF} is complete and agrees with the numerical experiments. However, the two-dimensional theory was based on a conjecture that has been proved only in a particular case (see Remark \ref{rem:2d_theorem} for more details) and therefore requires numerical validation.}

  In the present paper, we present focused numerical experiments to \as{validate} and
  extend the theoretical conclusions in \cite{BQCF,VKOr:blend2}.
\as{In particular, we study (i) whether stability of the B-QCF
method in 2D can be systematically improved with increasing the
blending width, (ii) whether a relatively narrow blending, as
suggested by the theory, is enough in practice, and (iii)
whether using the quintic spline (that has the regularity
assumed in the theory) has advantages over the cubic spline.}
  In addition we provide accuracy benchmarks similar to those in \cite{BQCEcomp}.
  Our numerical benchmarks demonstrate that the B-QCF scheme is a practical a/c coupling
  mechanism with performance (accuracy versus computational
  cost) superior to energy-based blending schemes.

\subsection{Summary}
In section~\ref{1DBQCFsection}, we introduce the B-QCF model for
a $1$D atomistic chain. We state the asymptotically optimal
condition on the blending size in Theorem~\ref{1DBQCFtheorem}
and apply a uniform expansion to the atomistic chain in subsection~\ref{1DoperaSubSec}.
The critical strain errors between the atomistic and B-QCF models with different blending
size are computed in this subsection. The numerical results perfectly match the analytic prediction,
that is, the errors decay polynomially in terms of the blending size.

In section~\ref{2DBQCFsection}, we establish the B-QCF model for
a $2$D hexagonal lattice. We state sufficient and necessary conditions
on the blending width under which the B-QCF operator is
positive definite. To \as{numerically} investigate the positive-definiteness of the B-QCF operators in $2$D,
we apply three different classes of deformations to the perfect lattice, which are the uniform
expansion, two types of shear deformation, and a general class of homogeneous deformations.
The results of $2$D uniform expansion are similar to those of the $1$D example, and they agree
with the theoretical conclusions well.

The stability regions of the different models
under homogeneous deformations are consistent with the analytic prediction.
By using a small blending region, the $2$D B-QCF operator becomes almost as stable as
the atomistic model, compared to the fact that
the stability region of the force-based quasicontinuum (QCF) method, i.e., the B-QCF method without blending region,
is a proper subset of the fully atomistic model \cite{qcf.iterative, qcf.stab, DobShapOrt:2011}.
However, the stability error under shear deformation for the B-QCF operator seems to only depend linearly on the system size,
which is observed from the numerical experiments.

In section~\ref{sec:accuracy}, we implement the B-QCF method
from a {\em practical} point of view. We briefly review the
accuracy results in terms of computational cost, i.e., the
total number of degrees of freedom  ${\rm DoF}$, and then
include some numerical experiments for a di-vacancy and
microcrack to demonstrate the superior accuracy of B-QCF over
other a/c coupling schemes that we have investigated previously
in \cite{BQCEcomp}.

\newpage

\section{The B-QCF Operator in $1$D.}\label{1DBQCFsection}
\subsection{Notation}
\label{notation}
We denote the scaled reference lattice by $\eps \mathbb{Z}:=
\{\eps\ell : \ell\in\mathbb{Z}\}$. We apply a macroscopic strain $F >
0$ to the lattice, which yields
\[
\mathbf{y}_F := F\eps \mathbb{Z} = (F \eps \ell)_{\ell \in \mathbb{Z}}.
\]
The space $\mathcal{U}$ of $2N$-periodic zero mean displacements
$\mathbf{u}=(u_{\ell})_{\ell \in \mathbb{Z}}$ from $\mathbf{y}_{F}$ is
given by
\[
\mathcal{U}:=\bigg\{\mathbf{u} : u_{\ell+2N}=u_{\ell}
\text{ for }\ell\in \mathbb{Z},
\text{ and }{\textstyle \sum_{\ell=-N+1}^{N}u_{\ell}}=0\bigg\},
\]
and we thus admit deformations $\mathbf{y}$ from the space
\[
\mathcal{Y}_{F}:=\{\mathbf{y}:
\mathbf{y}=\mathbf{y}_{F}+\mathbf{u}\text{ for some }\mathbf{u}\in
\mathcal{U}\}.
\]
We set $\eps=1/N$ throughout so that the reference length of the
computational cell remains fixed.

We define the discrete differentiation operator, $D\mathbf{u}$, on
periodic displacements by
\[
(D\mathbf{u})_{\ell}:=\frac{u_{\ell}-u_{\ell-1}}{\epsilon}, \quad
-\infty<\ell<\infty.
\]
We note that $\left(D\mathbf{u}\right)_{\ell}$ is also $2N$-periodic
in $\ell$ and satisfies the zero mean condition. We will often denote
$\left(D\mathbf{u}\right)_{\ell}$ by $Du_{\ell}$.
We then define $\left(D^{(2)}\mathbf{u}\right)_{\ell}$ and $\left(D^{(3)}\mathbf{u}\right)_{\ell}$
for $-\infty<\ell<\infty$ by
\begin{equation*}
\left(D^{(2)}\mathbf{u}\right)_{\ell}:=\frac{Du_{\ell+1}-Du_{\ell}}{\epsilon};\quad
\left(D^{(3)}\mathbf{u}\right)_{\ell}:={ \frac{D^{(2)}u_{\ell}-D^{(2)}u_{\ell-1}}{\epsilon}.}
\end{equation*}
To make the formulas more concise, we sometimes denote $Du_{\ell}$ by
$u'_{\ell}$, $D^{(2)}u_{\ell}$ by $u''_{\ell}$, etc., when there is no
confusion in the expressions.

For a displacement $\mathbf{u}\in \mathcal{U}$ and its discrete derivatives, we employ the weighted
discrete $\ell_{\epsilon}^{p}$  and $\ell_{\epsilon}^{\infty}$ norms by
{
\begin{align*}
\|\mathbf{u}\|_{\ell_{\epsilon}^{p}}&:= \left( \epsilon
\sum_{\ell=-N+1}^{N}|u_{\ell}|^{p}\right)^{1/p}\text{for }1\le p<\infty,
\qquad
\|\mathbf{u}\|_{\ell_{\epsilon}^{\infty}}:=\max\limits_{-N+1\le \ell\le N}|u_{\ell}|,
\end{align*}
}
and the weighted inner product { for $\ell_{\epsilon}^{2}$ is}
\[
\la \mathbf{u},\mathbf{w}\ra :=\sum\limits_{\ell=-N+1}^{N}\epsilon u_{\ell}w_{\ell}.
\]

\subsection{The B-QCF Operator}\label{1DoperaSubSec}
We consider a one-dimensional ($1$D) atomistic chain with  periodicity
$2N$, denoted ${\bf y} \in \mathcal{Y}_{F}$, under second-neighbor
pair interaction. The total atomistic energy
per period of ${\bf y}$ is given by
$\mathcal{E}^{a}(\mathbf{y})-\epsilon
\sum_{\ell=-N+1}^{N}f_{\ell}y_{\ell}$, where
\begin{equation}\label{AtomEnergy1D}
\mathcal{E}^{\rm{a}}(\mathbf{y})
=\epsilon\sum_{\ell=-N+1}^{N}\left[\phi(y'_{\ell})+\phi(y'_{\ell}+y'_{\ell-1})\right]
\end{equation}
{for external forces
  $f_{\ell}$ and a two-body potential $\phi \in C^2(0,
    +\infty)$ such \as{as
    the Morse potential 
    given by \eqref{MorsePotential}.}
    Implicitly we also assume that $\phi(r), \phi'(r)$ and
  $\phi''(r)$ decay rapidly as $r$
  increases, so that we only have to take into account first and second neighbors.}

The equilibrium equations are given by the force balance
at each atom: $F_\ell^a + f_\ell = 0$ where
\begin{align}\label{AtomEquil1D}
F_{\ell}^{\rm a}(\mathbf{y}):=\frac{-1}{\epsilon}\frac{\partial \mathcal{E}^{\rm a}(\mathbf{y})}{\partial y_{\ell}}
=& \frac{1}{\epsilon}\Big\{ \left[\phi'(y'_{\ell+1})+\phi'(y'_{\ell+2}+y'_{\ell+1})\right]
-\left[\phi'(y'_{\ell})+\phi'(y'_{\ell}+y'_{\ell-1})\right]
\Big\}.
\end{align}
{The linearized
equilibrium equations about $\mathbf{y}_{F}$ are
\[
\left(L^{\rm a}\mathbf{u}^{\rm a}\right)_{\ell}=f_{\ell},\quad\text{for}\quad\ell=-N+1,\dots,N,
\]
where $\left(L^{\rm{a}}\mathbf{v}\right)$ for a displacement $\mathbf{v}\in \mathcal{U}$ is given by
\[
\left(L^{\rm{a}}\mathbf{v}\right)_{\ell}:=\phi''_{F}\frac{\left(-v_{\ell+1}+2v_{\ell}-v_{\ell-1}\right)}{\epsilon^{2}}+
\phi''_{2F}\frac{\left(-v_{\ell+2}+2v_{\ell}-v_{\ell-2}\right)}{\epsilon^{2}}.
\]
Here and throughout
we use the notation $\phi''_{F}:=\phi''(F)$ and
$\phi''_{2F}:=\phi''(2F)$, where $\phi$ is the potential in
\eqref{AtomEnergy1D}.  We assume that $\phi''_{F} > 0$, which holds
for typical pair potentials such as the Lennard-Jones potential under
physically relevant deformations.
Appropriate extensions of the stability results in this paper can
likely be obtained for more general smooth deformations by utilizing
the more technical formalism developed, for example, in
\cite{bqce11,OrtnerShapeev:2010,ortner:qnl1d}.}

The local QC (or Cauchy-Born) approximation (QCL) uses the Cauchy-Born extrapolation rule \cite{Ortiz:1995a,Shimokawa:2004},
that is, approximating $y'_{\ell}+y'_{\ell-1}$  in \eqref{AtomEnergy1D} by $2y'_{\ell}$
in our context. Thus, the QCL energy is given by
\begin{equation}\label{QCLEnergy1D}
\mathcal{E}^{\rm{qcl}}(\mathbf{y})=\epsilon\sum_{\ell=-N+1}^{N}\left[\phi(y'_{\ell})+\phi(2y'_{\ell})\right].
\end{equation}
{
Then the local continuum forces $F^{\rm{qcl}}(\mathbf{y})$ are
\begin{align*}
F_{\ell}^{\rm{qcl}}(\mathbf{y}):=\frac{-1}{\epsilon}\frac{\partial \mathcal{E}^{\rm{qcl}}(\mathbf{y})}{\partial y_{\ell}}
=& \frac{1}{\epsilon}\Big\{ \left[\phi'(y'_{\ell+1})+2\phi'(2y'_{\ell+1})\right]
-\left[\phi'(y'_{\ell})+2\phi'(2y'_{\ell})\right]
\Big\}.
\end{align*}
}
We can similarly obtain the linearized QCL equilibrium equations about the uniform deformation
\[
\left(L^{\rm{qcl}}\mathbf{u}^{\rm{qcl}}\right)_{\ell}=f_{\ell}\quad \text{for}\quad \ell=-N+1,\dots, N,
\]
where the expression of $\left(L^{\rm{qcl}}\mathbf{v}\right)_{\ell}$ with $\mathbf{v}\in \mathcal{U}$
is
\[
\left(L^{\rm{qcl}}\mathbf{v}\right)_{\ell}:=
\left(\phi''_{F}+4\phi''_{2F}\right)\frac{\left(-v_{\ell+1}+2v_{\ell}-v_{\ell-1}\right)}{\epsilon^2}.
\]

The blended QCF (B-QCF) operator is obtained through smooth blending
of the atomistic and local QC models. Let $\beta : \mathbb{R} \to
\mathbb{R}$ be a ``smooth'' and $2$-periodic blending function, then
we define
\begin{displaymath}
  F_\ell^{\rm{bqcf}}(\mathbf{y}) := \beta_\ell F_\ell^{\rm a}(\mathbf{y}) +
  (1-\beta_\ell) F_\ell^{\rm{qcl}}(\mathbf{y}),
\end{displaymath}
where {$\beta_\ell := \beta(\epsilon \ell)$.} {Linearization} about
$\mathbf{y}_F$ yields the linearized B-QCF operator
\[
(L^{\rm{bqcf}}\mathbf{v})_{\ell}:=\beta_{\ell}
 (L^{\rm a}\mathbf{v})_{\ell}+(1-\beta_{\ell})(L^{\rm{qcl}}\mathbf{v})_{\ell}.
\]
Next, we define the blending \as{region $\mathcal{I}$ of width $K$:}
\begin{equation}\label{1DblendsizeDef}
\begin{split}
  \mathcal{I}:&=\big\{\ell \in {\helen{\{-N+1,\dots, N\} }} : 0 < \beta_{\ell+j} < 1
  \text{ for some } {\helen{j \in \{0, \pm 1, \pm 2\}}}
   \big\},\quad\text{ and }\\
  K:&= {\helen{\text{ the cardinality of the set }\mathcal{I}, }}
  \end{split}
\end{equation}
so that $D^{(j)}\beta_\ell = 0$ for all $\ell \in {\helen{\{-N+1,\dots, N\} }}
\setminus \mathcal{I}$ and $j \in \{1,2,3\}$. Thus $K$ is the size
of the compact support of $D^{(j)}\beta_\ell$. It is obvious that $K< 2N$.

\subsection{Positive-Definiteness of the B-QCF Operator}
We proved in \cite{BQCF} that the blending function
    $\beta$ can be chosen as a quintic polynomial such that
     \begin{itemize}
  \item[(i)] The $j$th derivatives of $\beta$ satisfy
    \begin{equation}
      \label{eq:BlendFunEst_upper}
      \| D^{(j)} \beta \|_{\ell^\infty} \leq C_\beta (K \varepsilon)^{-j},
      \quad \text{for } j = 1, 2, 3.
    \end{equation}
  \item[(ii)] This estimate is sharp in sense that, if $\beta_\ell$
    attains both the values $0$ and $1$, then
    \begin{equation}\label{BlendFunEst}
      \|D^{(j)}\beta\|_{\ell^{\infty}}\ge (K\varepsilon)^{-j},\quad
      \text{for } j=1,2,3.
    \end{equation}
 \end{itemize}

A linearized operator $L^{\rm w}$ with $\rm w \in \{{\rm a,c,bqcf}\}$,
is said to be positive definite in the $H^{1}$ norm or {\em coercive} if
there exists a constant $\gamma>0$ such that
\begin{equation}\label{PosDef}
\la L^{\rm w}\mathbf{u},\mathbf{u}\ra \ge\gamma
\|D\mathbf{u}\|_{\ell_{\varepsilon}^2}^2\quad \forall \mathbf{u}\in
\mathcal{U}.
\end{equation}
We have proved an asymptotically optimal stability condition on the blending region size
of the $1$D B-QCF operator in \cite{BQCF}.
\begin{theorem}\label{1DBQCFtheorem}
  Let $\mathcal{I}$ and $K$ be defined as in
  \eqref{1DblendsizeDef}, and suppose that $\beta$ is chosen to
  satisfy the upper bound \eqref{eq:BlendFunEst_upper}. Then there
  exists a constant $C_1 = C_1(C_\beta)$,
  such that
  \begin{equation}
    \label{eq:1d_coerc_lower}
    \langle L^{\rm{bqcf}} {\bf u}, {\bf u} \rangle \geq \big(c_0 - C_1
    |\phi_{2F}''| \big[ K^{-5/2} N^{1/2}\big]\big) \|D {\bf
      u}\|_{\ell^2_\epsilon}^2
    \qquad \forall {\bf u} \in \mathcal{U},
  \end{equation}
  where $c_0 = \min(\phi_F'', \phi_F'' + 4 \phi_{2F}'')$ is the
  atomistic stability constant.

  Moreover, if $\beta_\ell$ takes both the values $0$ and $1$, then
  there exist constants $C_2, C_3 > 0$, independent of $\mathcal{I}$,
  $N$, $\phi_F''$ and $\phi_{2F}''$, such that
  \begin{equation}
    \label{eq:1d_coerc_upper}
    \langle L^{\rm{bqcf}} {\bf u}, {\bf
      u} \rangle \leq \left(c_0 + \left\{C_2 - C_3 \left[ K^{-5/2} N^{1/2} \right]\right\}|\phi_{2F}''| \right)\|D {\bf
      u}\|_{\ell^2_\epsilon}^2\qquad \as{\text{for some } \mathbf{u}\in
\mathcal{U}\setminus\{0\}}.
  \end{equation}
\end{theorem}
From the conclusion of Theorem~\ref{1DBQCFtheorem}, we can immediately
get the following necessary and sufficient conditions on the blending
width $K$ for the operator $L^{\rm{bqcf}}$ to be coercive.
\begin{corollary}\label{1DBQCFcorollary}
Suppose that {$L^{\rm a}$} is positive-definite and that the blending function
{\helen{is sufficiently smooth}}. If the blending size $K$ satisfies $K\gg N^{1/5}$, then the B-QCF
operator $L^{\rm{bqcf}}$ is positive-definite and this estimate is
asymptotically optimal.
\end{corollary}
\subsection{$1$D Uniform Expansion Experiments.}

	We conduct numerical experiments in order to verify our theoretical findings.
	More precisely, we compare the decay rates of the error in critical strain as computed by B-QCF with the theoretically predicted rates as we increase the blending width $K$.

We use two kinds of blending functions: a cubic spline
\begin{equation}\label{GlobBlendFun}
\hat{B}\left(x\right)=\begin{cases} 0 &\quad x<0, \\
-2x^3+3x^2 &\quad 0\le x\le 1,\\
1&\quad x>1,
\end{cases}
\end{equation}
and a quintic spline
\begin{equation}\label{GlobQuinticFun}
\bar{B}(x)=\begin{cases} 0 &\quad x<0, \\
6x^5-15x^4+10x^3 &\quad 0\le x\le 1,\\
1&\quad x>1.
\end{cases}
\end{equation}
We scale $\hat{B}(x)$ and $\bar{B}(x)$ and define the blending functions for the atomistic chains as
\[
\hat{\beta}_{\ell}:=\hat{B}\left(\frac{\ell}{K}\right)\text{ and }
\bar{\beta}_{\ell}:=\bar{B}\left(\frac{\ell}{K}\right)\text{ for } \ell=-N+1,\dots,N. \]
Therefore, atoms with indices from $-N+1$ to $0$ belong to the
continuum region, from $1$ to $K-1$ belong to the blending region, and
from $K$ to $N$ belong to the atomistic region.
We note that $\bar{B}(x)$ has three bounded derivatives and hence it satisfies \eqref{eq:BlendFunEst_upper}, whereas for $\hat{B}(x)$ the second derivative has a jump, hence the third derivative does not exist.
Therefore, we expect that only $\bar{\beta}$ will yield the asymptotically
optimal stability estimates for the B-QCF method (see \cite{BQCF}).

For our interaction potential, we use the Morse potential
\begin{equation}\label{MorsePotential}
\phi(r)=\left[1-\exp(-\alpha(r-1))\right]^2,
\end{equation}
and we cut-off the interactions beyond the second nearest neighbor interactions.

We apply a uniform expansion to the atomistic chain: $\mathbf{y}_F := F\eps \mathbb{Z}$
with Dirichlet boundary condition:
\begin{equation}\label{Dirichlet1D}
u_{-N+1}=u_{N}=0.
\end{equation}
We then compute the critical strains of the atomistic and B-QCF models with
different blending size $K$ and fixed $N$.
The critical strains are defined as
\begin{equation}\label{1DcritDef}
\gamma^{\rm w}:=\max\left\{F>0: L^{\rm w}(\mathbf{y}_G)\text{ is
    positive definite for all } G \in [1,F)\right\},
\end{equation}
where $\rm{w}\in\{{\rm  a,c,bqcf}\}$ denotes the respective model.
\begin{remark}\label{rem:Dirichlet-to-periodic}
The stability bounds in Theorem~\ref{1DBQCFtheorem}
hold also for displacements $\mathbf{u}$ satisfying a homogeneous Dirichlet boundary condition.
To establish this, we note (1) that the bounds hold for constant
displacements as well, and (2) that any function satisfying \eqref{Dirichlet1D} can be extended to a periodic function (possibly with a nonzero mean).
Hence, Corollary~\ref{1DBQCFcorollary} also holds for
displacements $\mathbf{u}$ with homogeneous Dirichlet boundary conditions~\eqref{Dirichlet1D}.
\end{remark}

The computational results are shown in Figure~\ref{1DExpansionFig}.
In Figure~\ref{1DExpansionFig:quintic} we plot the dependence of the errors of quintic blending on $K$ for different values of $\alpha$.
We see that the graph of the error for the quintic blending is very close to the lower bound $K^{-5/2}$ as given by \eqref{eq:1d_coerc_lower} in Theorem~\ref{1DBQCFtheorem}.
Also, the error is lower for larger $\alpha$, which is also in accordance with the theoretical results.
Indeed, when $\alpha$ is large, the strength of the next-nearest neighbor interaction, $\phi''_{2F}$, is small relative to the nearest neighbor interaction $\phi''_{F}$, which contributes to a better stability of B-QCF according to \eqref{eq:1d_coerc_lower}.

Figure~\ref{1DExpansionFig:compare} shows the results of comparison of the cubic and the quintic blending.
We see that the cubic blending produces the error that seems to decay slower, like $K^{-2}$.
On the other hand, the quantitative difference between cubic and quintic is not large on the example considered.
To observe a significantly higher accuracy of the quintic spline, the computational domain size $N$ has to be much larger.
In addition, for larger $\alpha$, $N$ has to be even larger for the quintic blending to have advantage over the cubic blending.

\begin{figure}[h]
\begin{center}
\subfigure[Quintic spline blending\label{1DExpansionFig:quintic}]{\includegraphics[width= 8 cm]{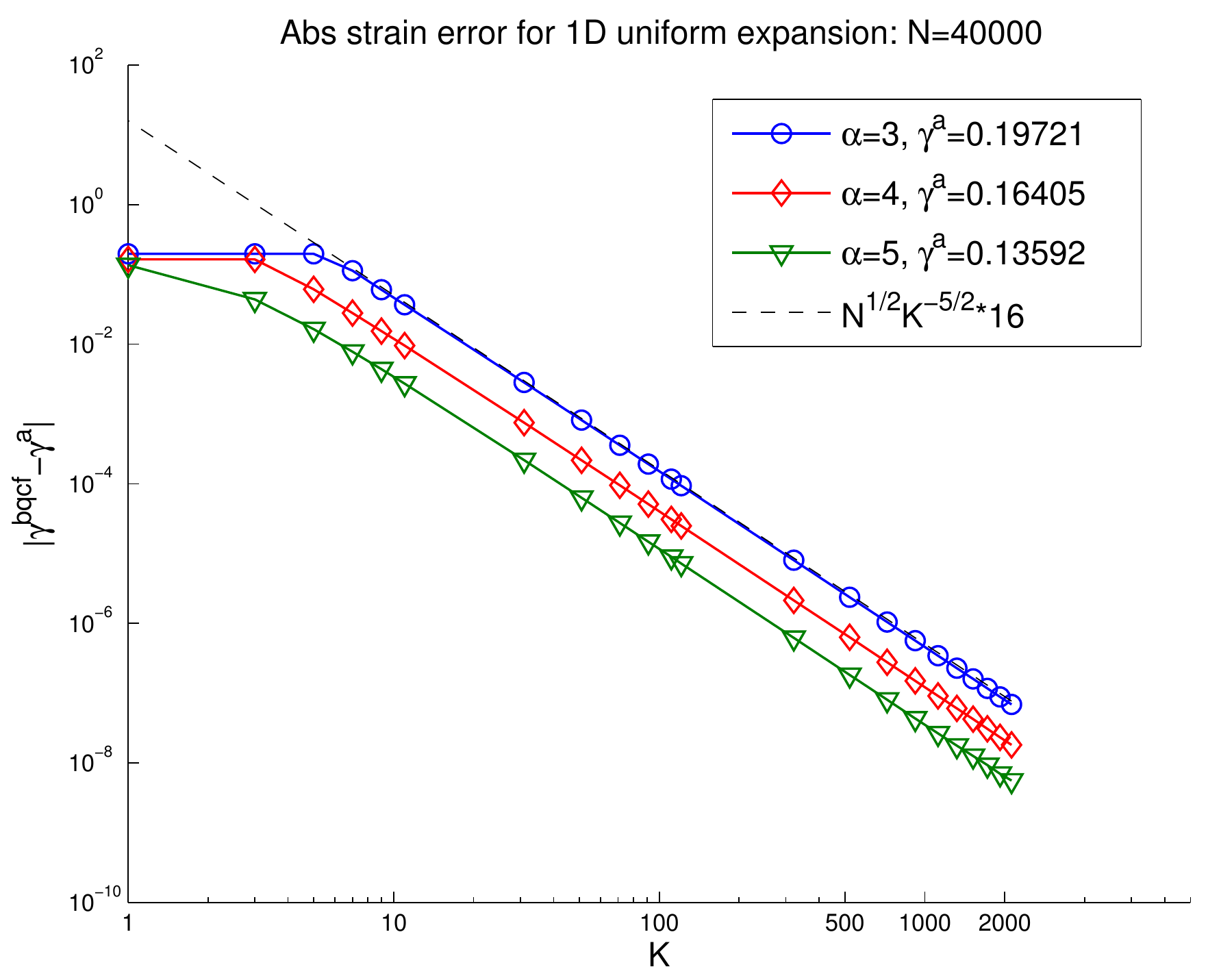}}
\subfigure[Quintic v.s. Cubic\label{1DExpansionFig:compare}]{\includegraphics[width= 8 cm]{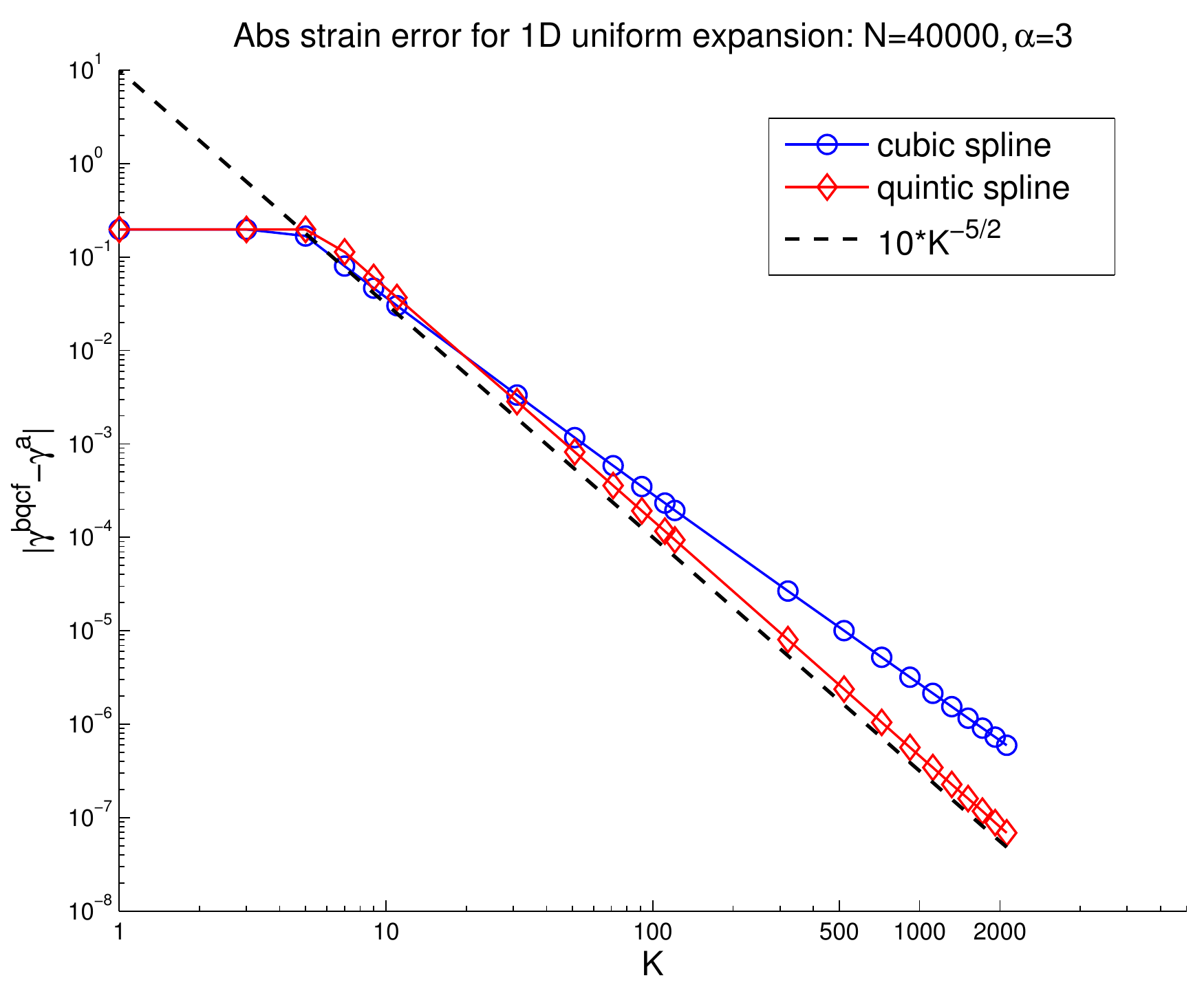}}
\end{center}
\caption{ (a) The absolute critical strain errors for a $1$D uniform
  expansion. We set $N=40,000,$  $\Delta \gamma=1/N^2$ where $\Delta\gamma$ is the strain increment used for
  testing stability, and $\gamma^{\rm{a}}$ and
  $\gamma^{\rm{bqcf}}$ are the critical strains for the atomistic and
  B-QCF models, respectively. The dashed line corresponds to the
  theoretical asymptote.  (b) The absolute critical strain errors of
  quintic and cubic blending functions with $N=40,000$ and $\alpha=3$.
  The solid line corresponds to the theoretical
  asymptote.}\label{1DExpansionFig}
\end{figure}


\section{The B-QCF Operator in $2$D.}\label{2DBQCFsection}

\subsection{The Triangular Lattice}\label{2DlatticeSubSec}
For some integer $N\in\mathbb{N}$ and $\epsilon:=1/N$, we define the
scaled 2D triangular lattice $\mathbb{L}$ to be
\[
\mathbb{L}:=\mathtt{A}_{6}\mathbb{Z}^2,\quad\text{where}\quad
\mathtt{A}_{6}:=\left[{a}_{1},{a}_{2}\right] :=
\epsilon\left[\begin{array}{cc}
1 & 1/2\\
0 & \sqrt{3}/2
\end{array}\right],
\]
where ${a}_{i},\,i=1,2$ are the scaled lattice vectors. Throughout
our analysis, we use the following definition of the periodic
reference cell
\[
\Omega:=\mathtt{A}_{6}(-N/2, N/2]^2\quad \text{and}\quad
\mathcal{L}:=\mathbb{L}\cap\Omega.
\]
We furthermore set ${a}_3=(-1/2\epsilon,
\sqrt{3}/2\epsilon)^{\mathtt{T}}$, then the set of
\emph{nearest-neighbor directions} is given by
\[
\mathcal{N}_{1}:=\{\pm{a}_1,\pm{a}_{2},\pm{a}_3\}.
\]
The set of \emph{next nearest-neighbor directions} is given by
\[
\mathcal{N}_{2}:=\{\pm {b}_1,\pm {b}_{2},\pm {b}_3\}, \quad
\text{where} \quad b_1:=a_1+a_2,\quad
b_2:=a_2+a_3,\quad\text{and}\quad b_3=a_3-a_1.
\]
We use the notation $\mathcal{N}:=\mathcal{N}_1\cup\mathcal{N}_2$ to
denote the directions of the neighboring bonds in the interaction
range of each atom (see Figure~\ref{AtomDomainFig}).

We identify all lattice functions $\mathbf{v} : \mathbb{L} \to
\mathbb{R}^2$ with their continuous, piecewise affine interpolants with
respect to the canonical triangulation $\mathcal{T}$ of $\mathbb{R}^2$
with nodes $\mathbb{L}$.

\begin{figure}[h]
\begin{center}
\subfigure[Neighbor set]{\includegraphics[height=5cm]{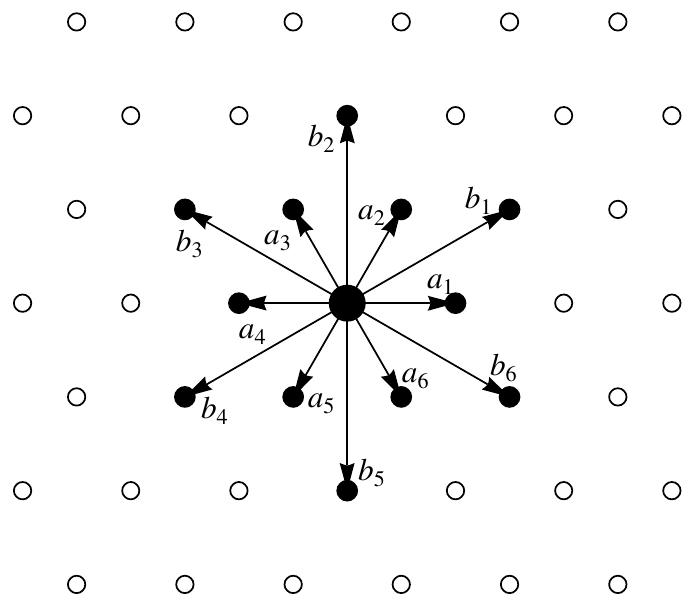}}\quad
\subfigure[Domain decomposition]{\includegraphics[height=5cm]{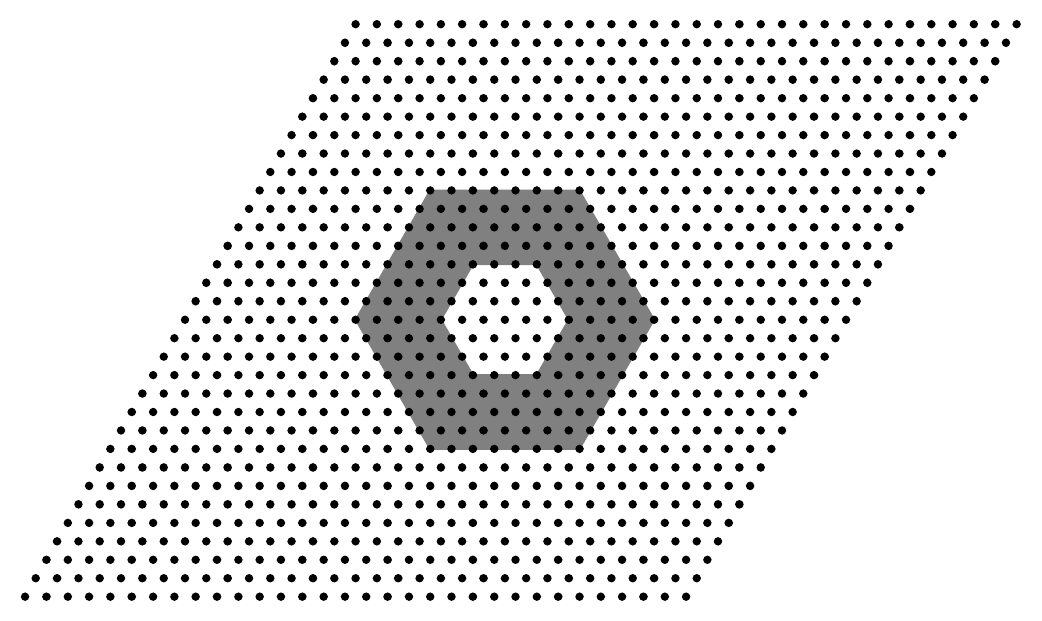}}
\end{center}
\caption{(a) The 12 neighboring bonds of each atom. (b) The periodic reference cell $\mathcal{L}:=\mathbb{L}\cap\Omega,$ the atomistic
  region $\Omega_a:=\mathtt{Hex}(\epsilon R_a),$ and the blending region
   $\Omega_b:=\mathtt{Hex}(\epsilon R_b)\setminus \Omega_a$. Here,
  $N=32,$ $R_a = 3$, $R_b = 7$, and $K = 4$.}\label{AtomDomainFig}
\end{figure}

\subsection{The Atomistic, Continuum, and Blending Regions}
Let $\mathtt{Hex}(R)$ denote the closed hexagon centered at the
origin, with sides aligned with the lattice directions $a_1, a_2,
a_3$, and diameter $2R$.  For $R_a< R_b\ <N  \in \mathbb{N}$, we define the
atomistic, blending, and continuum regions, respectively, as
\begin{displaymath}
 \Omega_a := \mathtt{Hex}(\epsilon R_a), \quad
  \Omega_b := \mathtt{Hex}(\epsilon R_b) \setminus \Omega_a,
  \quad \text{and} \quad
 \Omega_c:={\rm clos}\left(\Omega\setminus\left(\Omega_a\cup\Omega_b\right)\right).
\end{displaymath}
 We denote the blending width by $K := R_b - R_a$.
Moreover, we define the corresponding lattice sites
\begin{displaymath}
  \mathcal{L}^a := \mathcal{L}\cap \Omega_a, \qquad
  \mathcal{L}^b := \mathcal{L}\cap \Omega_b, \qquad \text{and} \qquad
  \mathcal{L}^c := \mathcal{L}\cap \Omega_c.
\end{displaymath}
For simplicity, we will again use $\mathcal{L}$ as the finite element
nodes, that is, every atom is a repatom.
For a map $\mathbf{v}:\mathbb{L}\rightarrow \mathbb{R}^2$ and bond
directions $r,s\in \mathcal{N}$, we define the finite difference
operators
\begin{displaymath}
D_{r}v(x):=\frac{v(x+r)-v(x)}{\epsilon}\quad\text{and}\quad
D_{r}D_{s}v(x):=\frac{D_{s}v(x+r)-D_{s}v(x)}{\epsilon}.
\end{displaymath}

We define the space of all admissible displacements, $\mathcal{U}$, as
all discrete functions $\mathbb{L}\rightarrow \mathbb{R}^2$ which are
$\Omega$-periodic and satisfy the mean zero condition on the computational domain:
\[
\mathcal{U}:=\Big\{\mathbf{u}:\mathbb{L}\rightarrow
\mathbb{R}^2 : \text{$u(x)$ is $\Omega$-periodic and
}{\textstyle \sum_{x\in\mathcal{L}}} u(x)=0 \Big\}.
\]
For a given matrix $B \in \mathbb{R}^{2 \times 2}$,
$\mathrm{det}(B)>0$, we admit deformations $\mathbf{y}$ from the space
\[
\mathcal{Y}_{B}:=\big\{ \mathbf{y}
:\mathbb{L}\rightarrow
\mathbb{R}^2: y(x)=Bx+u(x) {\helen{\,\,\, \forall x\in \mathbb{L},\,}}\text{ for some $\mathbf{u}\in\mathcal{U}$}
\big\}.
\]

For a displacement $\mathbf{u}\in \mathcal{U}$ and its discrete directional derivatives, we employ the weighted
discrete $\ell_{\epsilon}^{2}$  and $\ell_{\epsilon}^{\infty}$ norms given by
\begin{align*}
&\|\mathbf{u}\|_{\ell_{\epsilon}^{2}}:= \left( \epsilon^2
\sum_{x\in\mathcal{L}}|u(x)|^{2}\right)^{1/2},\quad
\|\mathbf{u}\|_{\ell_{\epsilon}^{\infty}}:=\max\limits_{x\in\mathcal{L}}|u(x)|,\quad\text{and}\\
&\qquad\quad\|D\mathbf{u}\|_{\ell_{\epsilon}^{2}}:=
\left(\epsilon^2\sum_{x\in\mathcal{L}}\sum_{i=1}^{3}|D_{a_i}u(x)|^2\right)^{1/2}.
\end{align*}
The inner product associated with $\ell^2_\epsilon$ is
\[
\la \mathbf{u},\mathbf{w}\ra :=\epsilon^2\sum\limits_{x\in\mathcal{L}}u(x)\cdot w(x).
\]

\subsection{The B-QCF operator.}
The total scaled atomistic energy for a periodic computational cell
$\Omega$ is
\begin{align}
\mathcal{E}^{\rm{a}}(\mathbf{y})=&\frac{\epsilon^2}{2}\sum_{x\in\mathcal{L}}\sum_{r\in \mathcal{N}}
\phi(D_{r}y({x}))\label{AtomEnergy2D}
= \epsilon^2 \sum_{x\in\mathcal{L}}\sum_{i=1}^{3}\big[\phi(D_{a_{i}}y({x}))
+\phi(D_{b_{i}}y({x}))\big],
\end{align}
where $\phi \in C^2(\mathbb{R}^2)$, for the sake of simplicity.
Typically, one assumes $\phi(r) = \varphi(|r|)$; the more general form
we use gives rise to a simplified notation; see also
\cite{OrtnerShapeev:2010}. We define $\phi'(r) \in \mathbb{R}^2$ and
$\phi''(r) \in \mathbb{R}^{2 \times 2}$ to be, respectively, the
gradient and hessian of $\phi$.

The equilibrium equations are given by the force balance at each atom,
\begin{equation}\label{AtomEquil2D}
F^{a}(x;y)+f(x;y)=0,\quad \text{for}\quad x\in\mathcal{L},
\end{equation}
where $f(x;y)$ are the external forces and $F^{a}(x;y)$ are the
atomistic forces (per unit area $\epsilon^2$)
\begin{align*}
F^{\rm a}(x;y):=&-\frac{1}{\epsilon^2}\frac{\partial \mathcal{E}^{\rm a}(\mathbf{y})}{\partial y(x)}\\
= &- \frac{1}{\epsilon} \sum_{i=1}^{3}\Big[\phi'\left(D_{a_i}y(x)\right) +\phi'\left(D_{-a_i}y(x)\right)
           \Big]
 -\frac{1}{\epsilon} \sum_{i=1}^{3}\Big[\phi'\left(D_{b_i}y(x)\right) +\phi'\left(D_{-b_i}y(x)\right)
           \Big].
\end{align*}
Again, since $\mathbf{u} = \mathbf{y} - \mathbf{y}_B$, where $y_B(x) =
B x$, is assumed to be small, we linearize the atomistic
equilibrium equation \eqref{AtomEquil2D} about $\mathbf{y}_B$:
\[
\left(L^{\rm a}\mathbf{u}^{\rm a}\right)(x)=f(x),\quad \text{for}\quad x\in \mathcal{L},
\]
where $\left(L^{\rm a}\mathbf{u}\right)(x)$, for a displacement $\mathbf{u}$, is given by
\[
\left(L^{\rm a}\mathbf{u}\right)(x)=-\sum_{i=1}^{3}\phi''(Ba_{i})D_{a_{i}}D_{a_{i}}u(x-a_{i})
-\sum_{i=1}^{3}\phi''(Bb_{i})D_{b_{i}}D_{b_{i}}u(x-b_{i}),\quad\text{for}\quad x\in\mathcal{L}.
\]

We use the Cauchy-Born extrapolation rule to approximate the nonlocal
atomistic model by a local continuum Cauchy-Born model
\cite{Ortiz:1995a,Shenoy:1999a,Miller:2003a}. Using the bond density
lemma \cite[Lemma 3.2]{OrtnerShapeev:2010} (see also
\cite{Shapeev2D:2011}), we can write the total QCL energy (the
discretized Cauchy-Born energy) as a sum of the bond density
integrals
\begin{equation}
  \label{QCLEnergy2D}
  \mathcal{E}^{\rm c}(\mathbf{y})= {\helen{\frac{1}{\Omega_0} } }\int_\Omega \sum_{r \in \mathcal{N}}
  \phi(\partial_r y) \, dx
  = \sum_{x\in \mathcal{L}}\sum_{r\in \mathcal{N}}
  \int_{0}^{1}\phi\big(\partial_{r}y(x+tr)\big)dt,
\end{equation}
where {\helen{the factor $\Omega_0:=\sqrt{3}/2$ is the volume of one primitive cell of $\mathcal{L}$}}
and $\partial_{r} y(x) := \frac{d}{dt} y(x+t r)|_{t = 0}$ denotes the
directional derivative. We compute the continuum force
\[F^{\rm c}(x;y) =
-\frac{1}{\epsilon^2} \frac{\partial\mathcal{E}^{\rm c}}{\partial y(x)},
\]
and linearize the force equation about the uniform deformation
$\mathbf{y}_{B}$ to obtain
\[
\left(L^{\rm c}\mathbf{u}^{\rm c}\right)(x)=f(x),\quad \text{for}\quad x\in \mathcal{L}.
\]

To formulate the B-QCF method, we let the blending function
$\beta(s):\mathbb{R}^2\rightarrow [0, 1]$ be a ``smooth'',
$\Omega$-periodic function.
Then, the (nonlinear) B-QCF forces are given through a convex combination of
$F^{\rm a}(x; y)$ and $F^{\rm c}(x; y)$:
\begin{displaymath}
  F^{\rm{bqcf}}(x; y) := \beta(x) F^{\rm a}(x; y) + (1-\beta(x)) F^{\rm c}(x; y),
\end{displaymath}
and linearizing the equilibrium equation $F^{\rm{bqcf}} + f = 0$ about
$y_B$ yields
\begin{equation}
  \label{eq:2}
  \begin{split}
    & (L^{\rm{bqcf}} \mathbf{u}^{\rm{bqcf}})(x) = f(x), \quad \text{for } x \in
    \mathcal{L},\\
    &\text{where} \quad  (L^{\rm{bqcf}} \mathbf{u})(x) = \beta(x) (L^{\rm a}
    \mathbf{u})(x) + (1-\beta(x)) (L^{\rm c}\mathbf{u})(x).
  \end{split}
\end{equation}

The $2$D blending function in our computational experiments will be defined radially using cubic and quintic splines:
\[
\hat{\beta}(x):=\hat{B}\left(\frac{\epsilon R_b-|x|}{\epsilon R_b-\epsilon R_a}\right)\quad\text{ and }\quad
\bar{\beta}(x):=\bar{B}\left(\frac{\epsilon R_b-|x|}{\epsilon R_b-\epsilon R_a}\right),
\]
where $\hat{B}(x)$ is given by \eqref{GlobBlendFun} and $\bar{B}(x)$ is given by \eqref{GlobQuinticFun}.
The function $\bar{\beta}(x)$ has the smoothness and satisfies 2D versions of the scaling bounds \eqref{eq:BlendFunEst_upper} needed for Theorem \ref{2DSuffTheorem} below, whereas $\hat{\beta}(x)$ does not
have a bounded third derivative.
We therefore can expect that $\hat{\beta}(x)$ will give a larger error asymptotically as compared to $\bar{\beta}(x)$.

\subsection{Positivity of the B-QCF operator in 2D}
Necessary and sufficient conditions for $L^{\rm{bqcf}}$ to be
positive-definite are given in \cite{BQCF}. To make this paper more
concise, we only state the conclusions without proof. First, we state
{\helen{a lower bound}} for $\la L^{\rm{bqcf}} \mathbf{u}, \mathbf{u} \ra$:

\begin{theorem}\label{2DBlendSizeThm}
  Suppose that $\beta\in C^3$
  and satisfies the scaling bounds \eqref{eq:BlendFunEst_upper}; then,
  \begin{displaymath}
    \la L^{\rm{bqcf}} \mathbf{u}, \mathbf{u} \ra \geq \gamma_{\rm{bqcf}}
    \| D\mathbf{u} \|_{\ell^2_\epsilon}^2,
  \end{displaymath}
  where
  \begin{equation}\label{asympError}
    \gamma_{\rm{bqcf}} := \tilde{\gamma} -
    C\, \big[ K^{-5/2} R_b^{1/2} |\log(R_b/N)|^{1/2} \big],
  \end{equation}
  where $C$ is a generic constant independent of $N$,
  and $\tilde{\gamma}$ is the coercivity constant for
  the operator $\tilde{L}$:
  \begin{equation*}
  \la \tilde{L} \mathbf{u}, \mathbf{u} \ra :=
  \la L^{\rm c} \mathbf{u}, \mathbf{u} \ra - \epsilon^4 \sum_{i = 1}^3
  \sum_{x \in \mathcal{L}} \beta(x-a_2) \big| D_{a_i}
  D_{a_{i+1}} u(x - a_1 - a_2) \big|_{b_i}^2  \geq \tilde{\gamma} \| D \mathbf{u}
    \|_{\ell^2_\epsilon}^2 \qquad \forall \mathbf{u} \in \mathcal{U}.
\end{equation*}
\end{theorem}
One can see very clearly that, whenever $N$ is polynomial in $R_b$
  and $K \gg R_b^{1/5}$, then $L^{\rm bqcf}$ can be expected to be
  coercive. Both are natural and easy to achieve.  We can thus deduce
the following result for the coercivity of $L^{\rm{bqcf}}:$

\begin{corollary}\label{2DSuffTheorem}
  Suppose that $\tilde{L}$ is positive-definite and that the
  blending function
      $\beta\in C^3$ and satisfies the scaling bounds \eqref{eq:BlendFunEst_upper}. Let the number of atoms $R_a$ along the
    radius be of order $N^{\alpha}$ with $0\le \alpha \le 1$. If the
  blending width $K$ satisfies
\begin{align}
K\gg
\begin{cases}
|\log{N}|^{1/4}, \quad \alpha=0,\nonumber\\
|\log{N}|^{1/5}N^{\alpha/5},\quad 0<\alpha<1,\label{2DsuffCond}\\
N^{1/5},\quad \alpha=1,\nonumber
\end{cases}
\end{align}
then the B-QCF operator $L^{\rm{bqcf}}$ is positive-definite.
\end{corollary}
\begin{remark}\label{rem:2d_theorem}
$\mathstrut$
\begin{itemize}
\item[(a)]
  The stability result of Theorem \ref{2DBlendSizeThm}, and hence of
  Corollary \ref{2DSuffTheorem}, is \as{based on the conjecture that the
  operator $\tilde{L}$ is stable}. In \cite{BQCF} we show that $\tilde{L}$ is
  indeed stable whenever nearest-neighbor interactions dominate.

  Moreover, based on the analysis and numerical experiments in
  \cite{OrtnerShapeev:2010} for a similar linearized operator, we
  expect that the region of stability for $\tilde{L}$ is the same as
  for $L^a$ as $N,\, R_a,\, R_b \rightarrow \infty$. We therefore expect
  that the result of Theorem~\ref{2DSuffTheorem} holds (up to a
  controllable error) if coercivity of $\tilde{L}$ is replaced by
  coercivity of $L^{\rm{a}}$ in the hypothesis.

\item[(b)] One can apply the argument of Remark
  \ref{rem:Dirichlet-to-periodic} to conclude that the results of
  Theorem \ref{2DBlendSizeThm} and Corollary \ref{2DSuffTheorem} are
  valid for homogeneous Dirichlet boundary conditions as well.
\end{itemize}

\end{remark}

By constructing a radial counterexample similar to our $1$D
counterexample, we can observe that our conditions in Corollary
\ref{2DSuffTheorem} are essentially necessary.

\begin{theorem}\label{2DNeccTheorem}
Suppose that $L^{\rm a}$ is positive-definite and that
the blending function $\beta\in C^3$
and satisfies the scaling bounds \eqref{eq:BlendFunEst_upper}.
The number of atoms $R_a$ along the radius is of order
$N^{\alpha}$ with $0< \alpha \le 1$. If the blending width $K$ is
$K\ll N^{\alpha/5}$, then the B-QCF operator
$L^{\rm{bqcf}}$ {cannot be} positive-definite and we can construct a radial counterexample in this case.
\end{theorem}
We note that there is a gap between the necessary and sufficient conditions for $0<\alpha<1$.
In addition, we have no necessary condition for $\alpha=0$, which corresponds to
a fixed atomistic core independent of
the reference cell $\Omega$.
\subsection{2D numerical experiments for B-QCF operators.}\label{2DNumerSection}
In this subsection, we will continue the numerical experiments for the $2$D B-QCF models to verify the theoretical
findings by comparing the decay rates of the error
in critical strain as computed by B-QCF with
the theoretically predicted rates as we increase the blending width $K$.

\begin{enumerate}
\item{\bf Uniform expansion.}

We first consider the simplest $2$D deformation: we apply a uniform expansion $y(x)=Bx$ with
\[
B=\gamma\left(
\begin{array}{cc}
1\quad&0\\
0\quad&1
\end{array}
\right)
\] to the perfect lattice $\mathcal{L}$ with Dirichlet boundary
condition:
\begin{equation}\label{Dirichlet2D}
u(x)=0\quad \forall x\in \partial\Omega.
\end{equation}
Then we compute the critical strains $\gamma$ of
the atomistic and B-QCF models with different blending region
width $K$.

We note that the $2$D conclusions also depend on the size
of the atomistic region. Therefore we let $R_a=K^{5/3}$ in order to narrow the
dependence only to the blending width $K$. Then the asymptotical term in \eqref{asympError}
for sufficiently large $N$ is approximately
{\helen{
\begin{align*}
 K^{-5/2} R_b^{1/2}| \log(R_b/N)|^{1/2}
= & K^{-5/2} (R_a+K)^{1/2} | \log(R_b/N)|^{1/2} \\
\approx &  K^{-5/2} R_a^{1/2}=K^{-5/3}=R_a^{-1},
\end{align*}
}}
which means that the error in $\gamma_{\rm{bqcf}}$ is systematically
improvable.

The choice of scaling $R_a=K^{5/3}$ is motivated by the results in
  \cite{HudsOrt:a} which indicate that, generically, one should expect
  an $O(R_a^{-1})$ error in the regions of stability between the
  infinite lattice atomistic model and the atomistic model in a domain
  with radius $R_a$. {\helen{In the computation, we assign integer values for $K$
  and use the rounded values for $R_a$, that is $R_a=\lfloor K^{5/3} \rfloor$.}}

The critical strains are defined as
\begin{equation}\label{2DcritDef}
\gamma^{\rm w}:=\max\big\{\bar\gamma>0: L^{\rm w}(B\mathbf{x})\text{
    is positive definite for } \gamma \in [0,\bar\gamma)\big\},
\end{equation}
where $\rm{w}\in\{\rm a,\rm{bqcf}\}$ denote the models. Here we use the MATLAB function {\tt eigs} \cite{MATLAB:2012} to compute the
smallest eigenvalue of the
symmetric part of $L^{\rm w}(B\mathbf{x})$ and thus determine the
 positive-definiteness of $L^{\rm w}(B\mathbf{x})$.

We also define the increment of the strain $\gamma$ in each step by $\Delta \gamma$.
The results in \cite{OrtnerShapeev:2010,HudsOrt:a,DobShapOrt:2011,doblusort:qce.stab}
suggest that the theoretical increments be of order $O(N^{-2})$ (at least, for finding the critical strain of a uniform lattice), and we set $\Delta \gamma=10^{-8}$
which is sufficiently small considering $N=200$ or $300$ in our experiments.

\begin{figure}[h]
\begin{center}
\subfigure[Quintic spline blending]{\includegraphics[width= 8 cm]{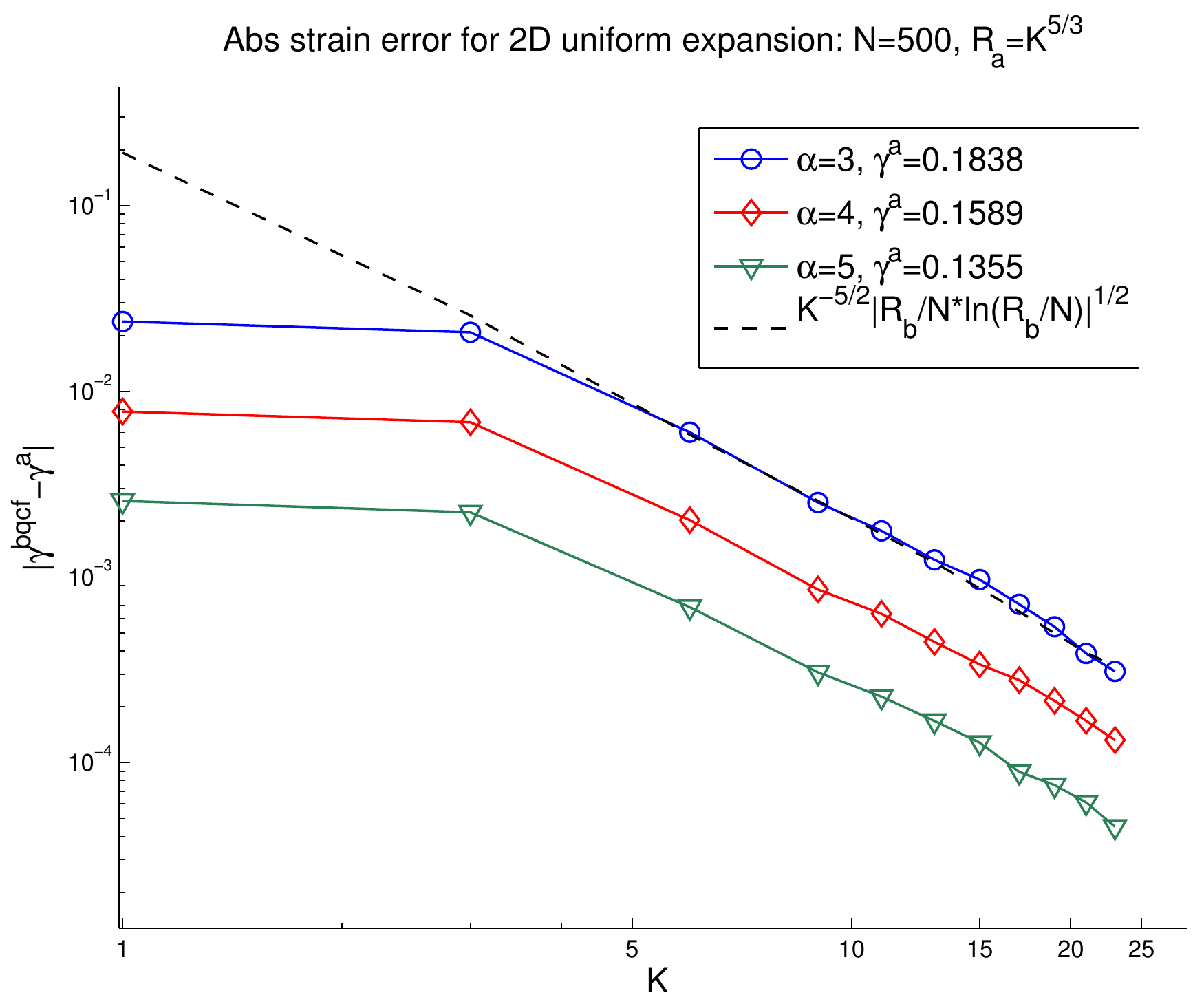}}
\subfigure[Quintic v.s. Cubic\label{2DExpanFigb}]{\includegraphics[width= 8 cm]{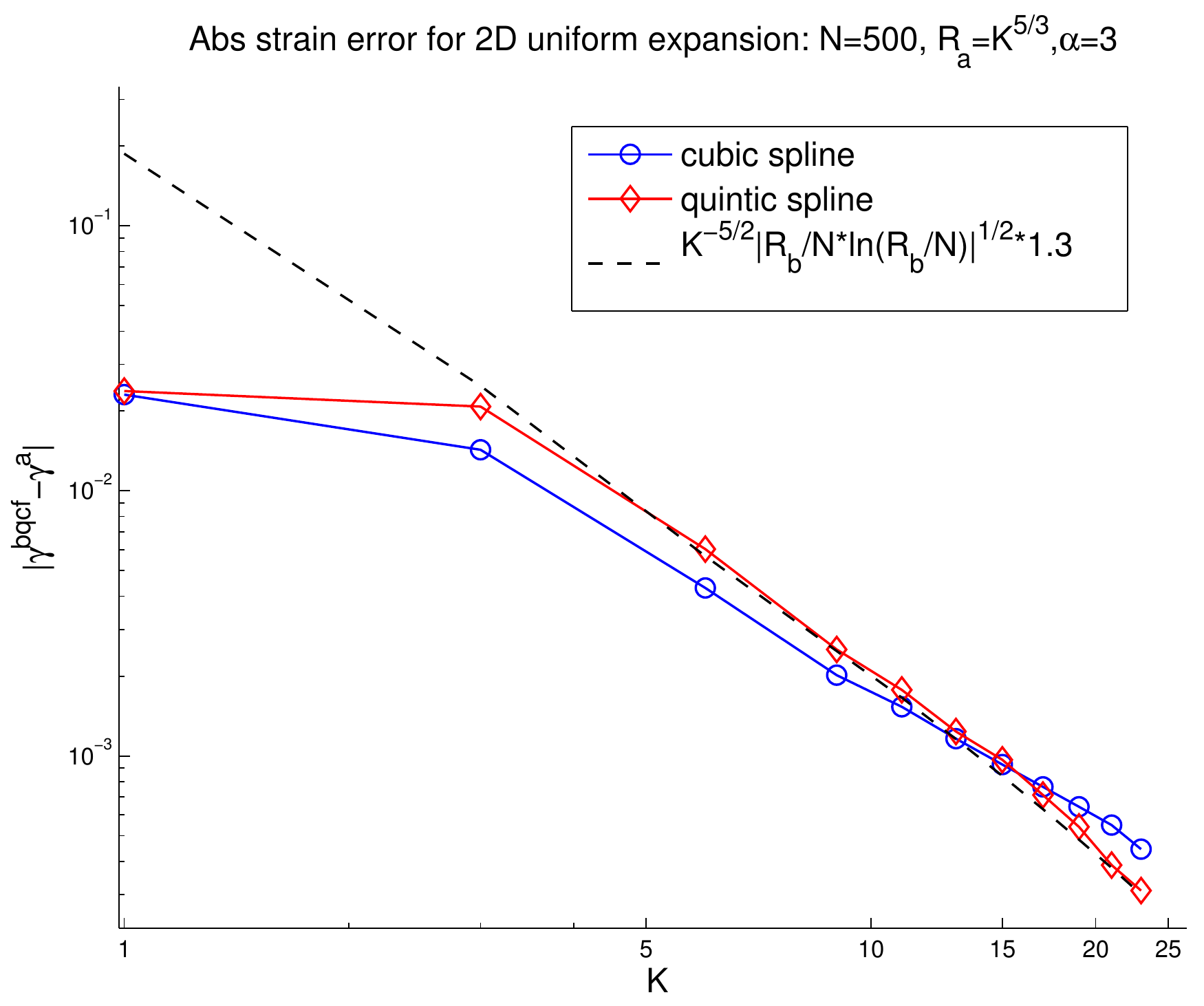}}
\end{center}
\caption{(a) The absolute critical strain errors for the $2$D uniform
expansion. {\helen{We set $N=500$}}, and we denote
the critical strains for the
atomistic and B-QCF models by
$\gamma^{\rm a}$ and $\gamma^{\rm{bqcf}},$
respectively. The dashed line corresponds to the theoretical asymptote. (b) The absolute critical strain errors for the quintic and cubic blending functions
with {\helen{$N=500$}} and $\alpha=3$. The solid line corresponds to the theoretical asymptote.}\label{2DExpanFig}
\end{figure}
We plot the difference of the critical strains with different blending width $K$ in Figure~\ref{2DExpanFig}.
The numerical critical strain
errors in the left figure approach the analytical asymptote as $K$ increases.
There are larger fluctuations of errors as compared to the $1$D case,
which is likely due to round-off errors in calculating {\helen{$R_a=K^{5/3}$}}.
Thus, the slopes of the errors with quintic blending agree
with the theoretical prediction in Theorem~\ref{2DBlendSizeThm}.
Also, similarly to the 1D results, the error is smaller
when the nearest neighbor interaction {\helen{dominates}} (that is, when $\alpha$ is large).

Although
the slope of the errors with cubic blending seems to be one half order less than that with quintic blending (see Figure \ref{2DExpanFigb}),
the computed errors for cubic blending are slightly smaller for the {\helen{relatively small $N$ considered.}}
We expect that for a sufficiently large system, the quintic blending would be more accurate.
In addition, the $2$D errors for uniform expansion  are similar to the $1$D results.
This is reasonable since the $2$D uniform expansion is similar to the $1$D deformation.

\item{\bf Uniform shear deformation.}

  We now investigate stability of B-QCF under shear deformation. We
  apply a y-directional shear deformation to the hexagonal lattice
  $\Omega$ with Dirichlet boundary conditions \eqref{Dirichlet2D}. The
  y-directional shear is $y(x)=\tilde{B}x$ with
\[
\tilde{B}=\left(
\begin{array}{cc}
1\quad&0\\
\gamma\quad&1
\end{array}
\right).
\]
The critical strain errors between the B-QCF and atomistic models with the quintic blending
are plotted in Figure~\ref{2DShearFig}.
\begin{figure}[h]
\begin{center}
\subfigure[ Error vs $N$]{\includegraphics[width= 7.7 cm]{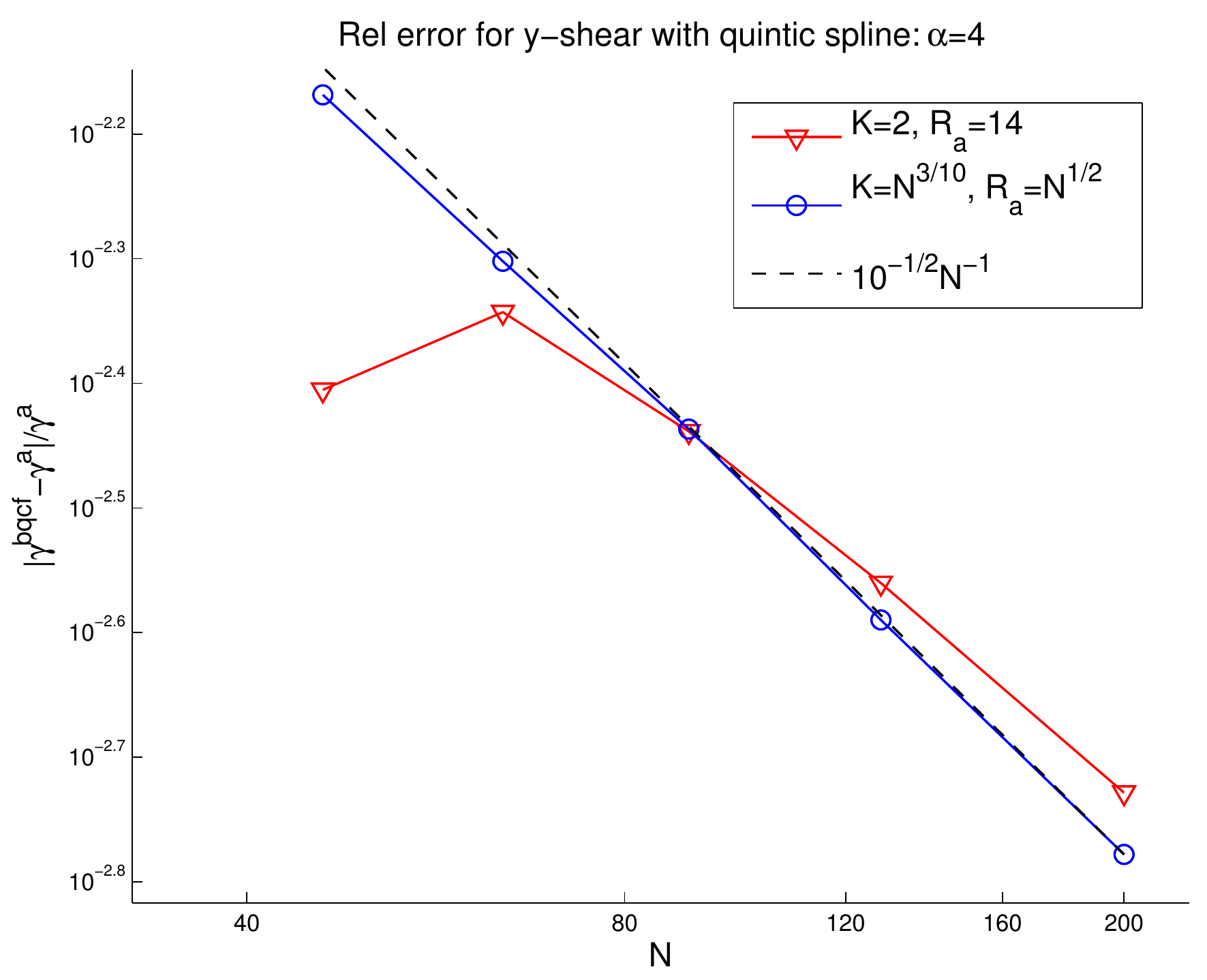}}\vspace{-2pt}
\subfigure[ Error vs $R_a$]{\includegraphics[width= 7.7 cm]{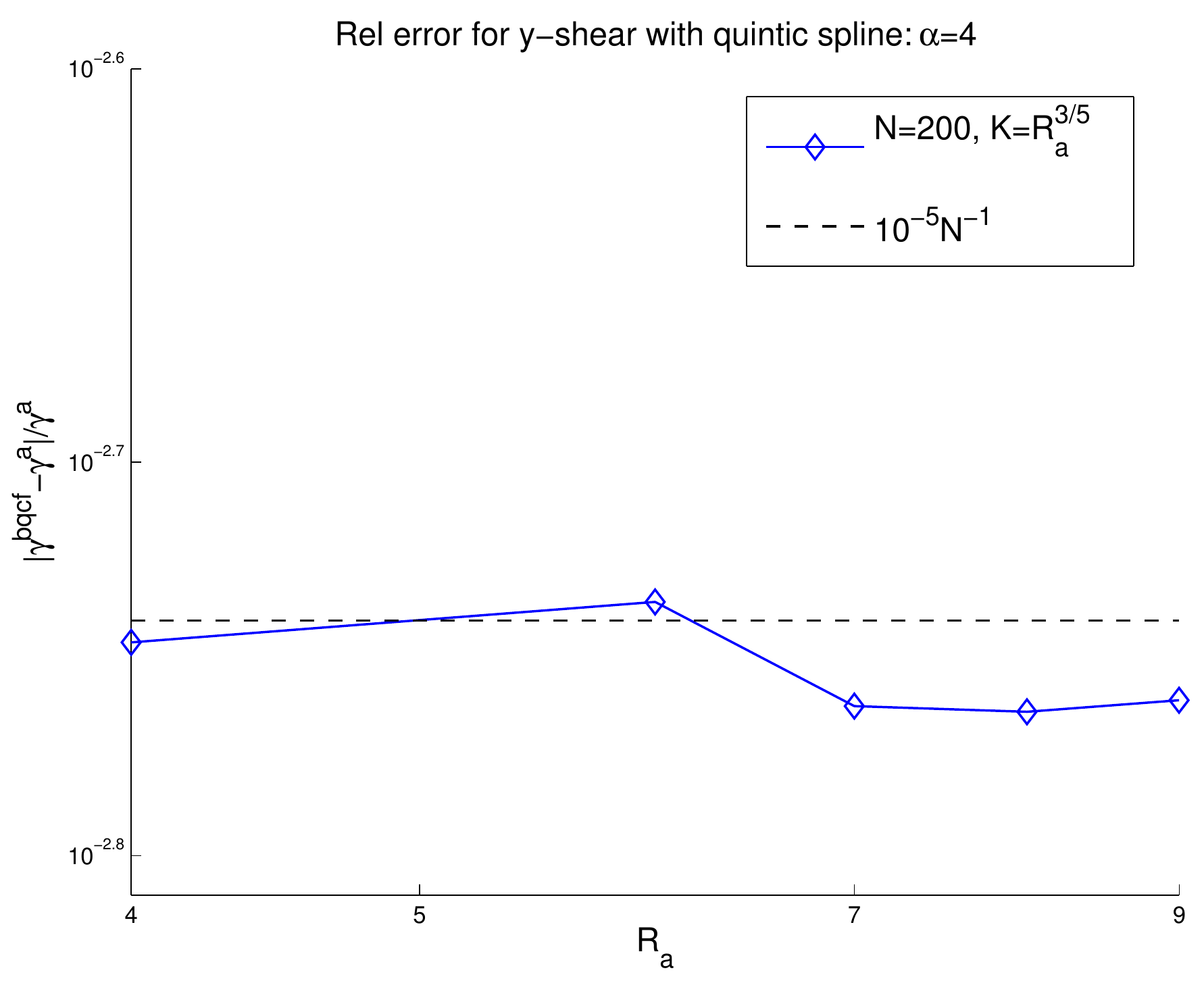}}
\end{center}
\caption{The relative critical strain error for the y-directional shear
deformation.
$\gamma^{\rm a}$ and $\gamma^{\rm{bqcf}}$ are the critical
strains for the atomistic and B-QCF models
respectively. {\helen{For $N=200$, $\gamma^a \approx 0.1813.$}}
The dashed line corresponds to the theoretical asymptote.
The fluctuations in the plotted error for $N=\const$ seems to be due to round-off errors in calculating $R_a$ and $K$.
{\helen{The method parameters were rounded as follows: in (a) $K=\lfloor N^{3/10}\rfloor$ and $R_a=\lfloor N^{1/2}\rfloor$, and in (b)
$K=\lfloor R_a^{3/5}\rfloor$.}}
 }\label{2DShearFig}
\end{figure}

In Figure~\ref{2DShearFig} we plot the critical strain
errors in the following three regimes: (1) $N$ increases,
$R_a=\const$, $K=\const$, (2) all three parameters
increase, and (3) $N=\const$, $R_a$ and $K$ increases.
\as{The choice of constant parameters, $K=2$ and $R_a=14$,
does not follow the scaling $K \approx R_a^{3/5}$, and was
made to show that such nonoptimal parameters do not
significantly affect the results in this case.} The results
indicate that the error in this case depends on $N$, but
does not depend on $R_a$ or $K$. This means that, for shear
deformations, the local continuum approximation and its
finite element coarse-graining contributes most of the
error.

We explain such a qualitative difference between the uniform expansion
and the shear deformation in the following way.  For the uniform
expansion the onset of instability is due to competition of
interaction of the nearest neighbors (NNs), contributing to stability,
and the second nearest neighbors (NNNs), contributing to instability.
On the other hand, for the shear deformation the onset of
  instability is primarily due to competition between elongated and
  compressed NN bonds.  Therefore, for the uniform expansion it is
  important to reduce the interface error which distorts the NNN
  interaction, whereas in shear deformation the NNN interactions do
  not contribute significantly to stability errors. Since, for NN
  interaction, the atomistic, Cauchy-Born and B-QCF models are
  identical, the stability error only depends on the domain size.

\item{\bf Regions of stability.}

We now combine the uniform expansion and shear deformation together and
study the stability region of $L^{\rm{bqcf}}$ for a general
class of homogeneous deformations.
We consider the following \as{family} of deformations which
involve shear, expansion, and compression.
\[
B=\left(\begin{array}{cc}
1+s\quad& 0.1\\
0 \quad&1+r
\end{array}
\right).
\]
Applying these specific homogeneous deformations to the hexagonal
lattice in the reference cell and again using the Dirichlet boundary
condition, we plot the stability regions (regions where the operators
are positive definite) in Figure~\ref{StabRegRa50Fig1}.
\begin{figure}[h]
\begin{center}
\subfigure[$\alpha=2$]{\includegraphics[height=7.1cm]{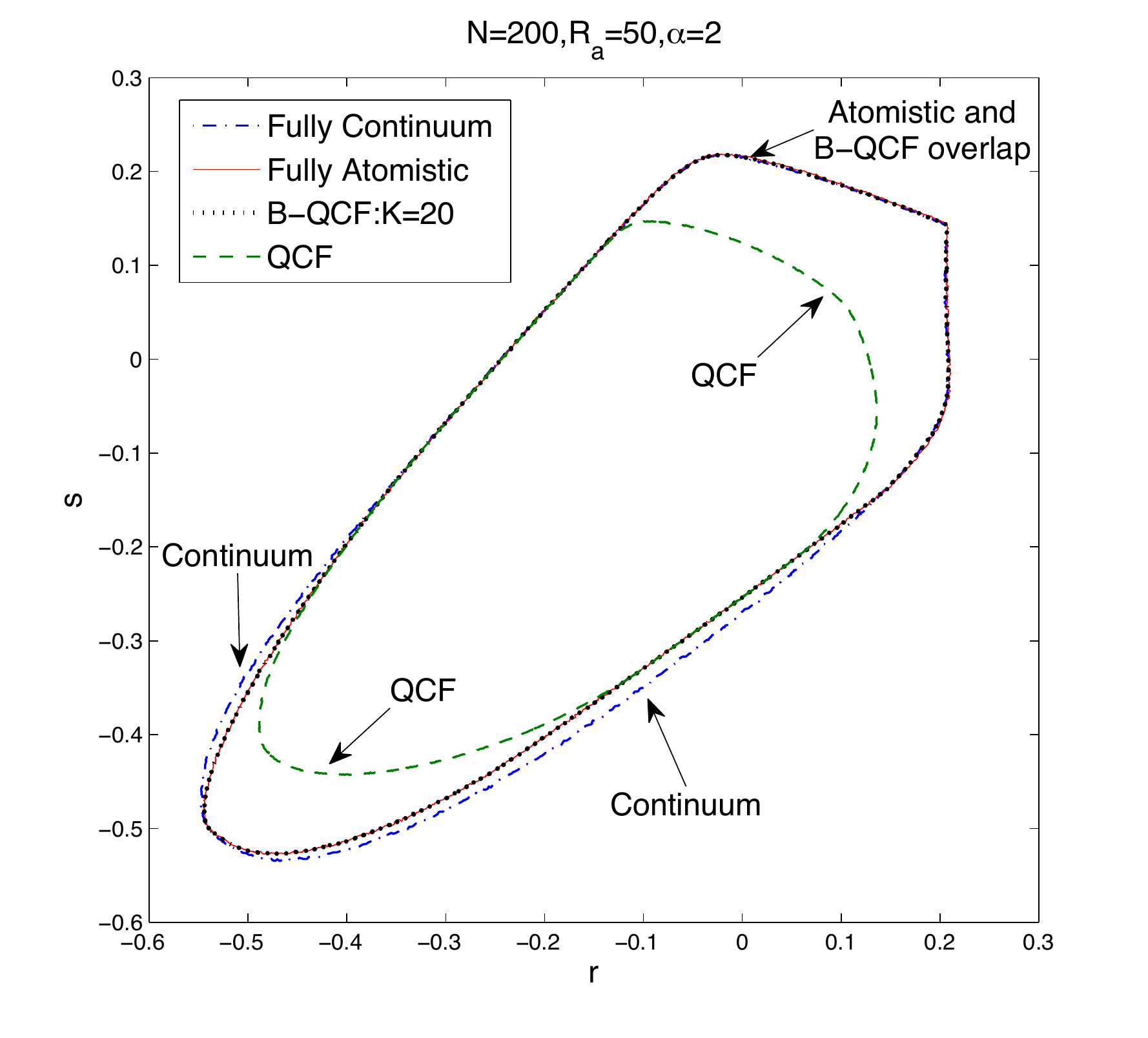}}
 \subfigure[$\alpha=4$]{\includegraphics [height =7.1cm]{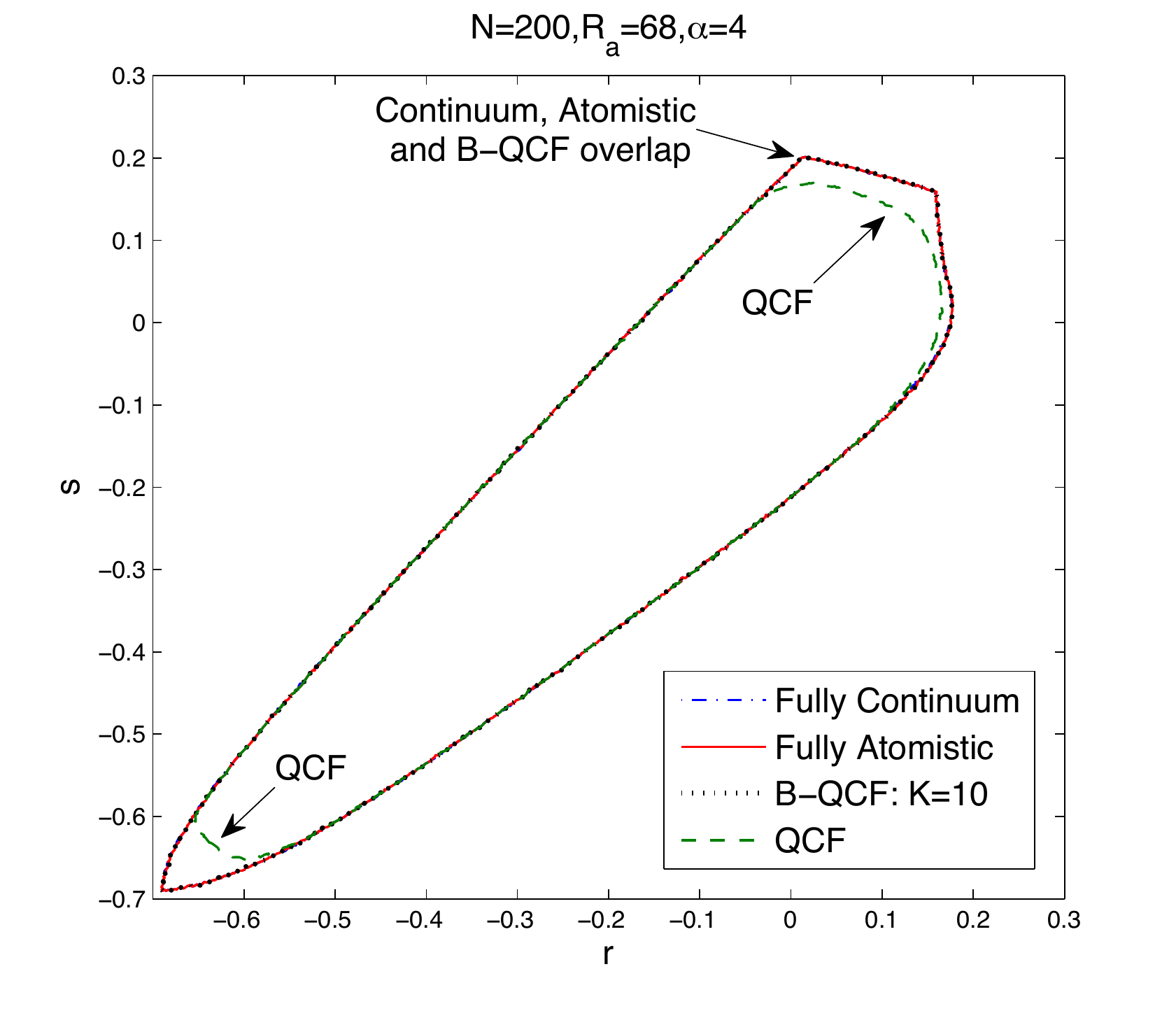}}
\end{center}
\caption{The stability regions of the different models.
These closed curves are the boundaries of the stability regions for the atomistic, B-QCF, and the local continuum models,
respectively. The curves with indicators are for QCF.}\label{StabRegRa50Fig1.pdf}
\end{figure}

We observe that the stability regions of the B-QCF model with
different blending sizes are all proper subsets of the atomistic
model.  In addition, the fully atomistic and continuum models are very
close to each other, which agrees with the stability analysis of the
perfect lattice \cite{HudsOrt:a}.  Also, when $\alpha$ increases,
which means the next-nearest neighbor interactions become less
important, the difference becomes smaller.

There is a visible difference in the stability regions between the QCF
model and the exact atomistic model, whereas the difference between
the B-QCF model and the atomistic model is almost not seen.  This
implies that using a blending region can significantly improve the
stability properties of the approximation models.

\item{\bf Stability of micro-cracks.}

The experiments that we have reported up to this point were based on perfect lattices. Now we
apply the B-QCF model to lattices with local defects.

The atomistic system is as follows.  There is a micro-crack in the
center of the domain $\Omega$ with length $5$, i.e., $5$ atoms are
removed from the lattice (see Figure~\ref{Fig:crackEq}).
\as{Hence, we redefine accordingly the positions of atoms in the reference configuration $\mathbf{x}$, the interaction energy $\mathcal{E}^{a}$, etc.}
We impose a
vertical stretching $B= \left[\begin{array}{cc}
    1&0\\
    0&1+\gamma
\end{array}
\right]$
on the lattice and compute the critical strains $\gamma^\c>0$ beyond which the system loses stability.

\begin{figure}[h]
\begin{center}
\includegraphics[width=16 cm ]{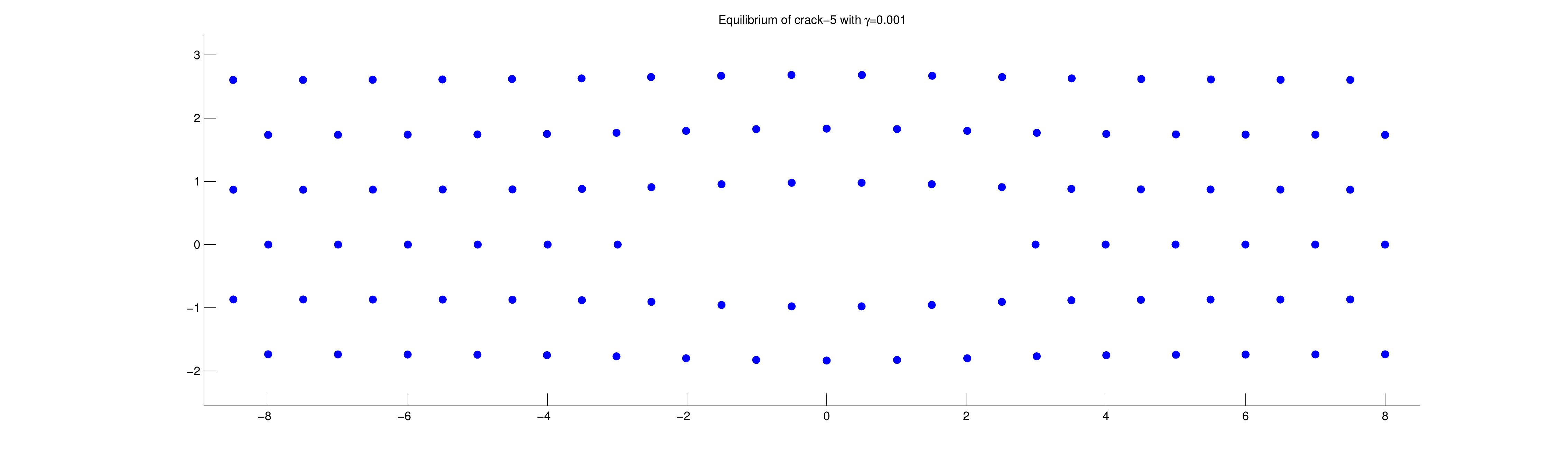}
\end{center}
\caption{The stable equilibrium configuration of the micro-crack with crack length$=5$ and
$\gamma=0.001$, and the $\ell^{\infty}_{\epsilon}$ norm of
the force residual is of order $O(10^{-12}).$ }\label{Fig:crackEq}
\end{figure}
We computed the critical strain $\gamma^\c$ in the following way.
Given $\gamma>0$, we use Newton's iteration method to solve the following force equations for $\mathbf{y}^{\gamma}$
with the initial guess $\mathbf{y}_F=B\mathbf{x}$:
\[
F^{w}(x;y^{\gamma})=0\quad\text{for}\quad x\in \Omega\setminus\partial \Omega.
\]
We set the tolerance for the $\ell^{\infty}_\epsilon$ norm of the force residual
of the Newton's iteration to be $10^{-5}$.
To prevent the configuration from ``jumping out'' of the local energy well corresponding to the defect under consideration,
we require at each step the $\ell_{\epsilon}^{\infty}$ norm of force residual
to be less than $100$ and the positive-definiteness of $L^{w}(\mathbf{y})$,
where $\mathbf{y}$ is the current configuration.
If any of the two requirements is not met, then the current $\gamma$ is regarded as an unstable strain.
When the force residual is smaller than the tolerance, the configuration $\mathbf{y}^{*}$ is thought to be in its equilibrium
of the local energy well. Then we check the positive definiteness of corresponding operator $L^{\rm{w}}(\mathbf{y}^*)$ with
the equilibrium configuration $\mathbf{y}^*$.
The nonlinear critical strain is thus defined as
\begin{equation*}
\gamma^\c:=\max\{\bar\gamma>0:L^{\rm{w}}(\mathbf{y}^*)\text{ is
  positive definite for }\gamma \in [0,\,\bar\gamma)\}.
\end{equation*}

The plot of critical strain for the B-QCF models are shown in Figure~\ref{Fig:Crack}.
\begin{figure}[h]
\begin{center}
\includegraphics [height =7.5 cm]{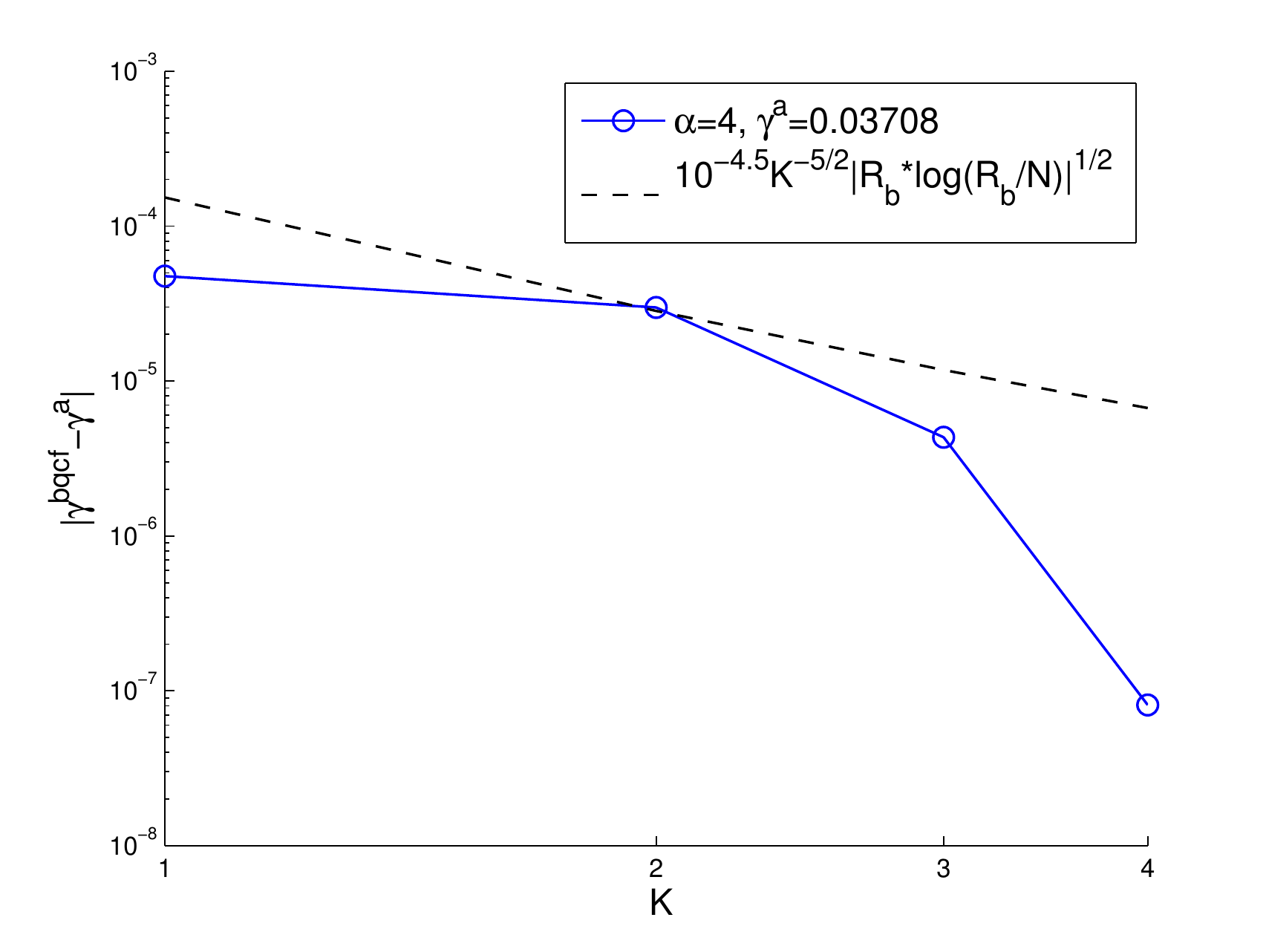}
\end{center}
\caption{The nonlinear critical strain error
for vertical stretching. We set $N=200$, crack length$=5$, and $R_a=\max\{K^2, 6\}$.
$\gamma^{\rm{a}}$, $\gamma^{\rm{bqcf}}$  are the critical
strains for the atomistic and B-QCF models,
respectively. The dashed line corresponds to the theoretical asymptote.}\label{Fig:Crack}
\end{figure}
\as{Even though we choose the blending width $K\approx R^{1/2}$ slightly smaller then our previous choice ($K\approx R^{3/5}$), }%
we observe the nonlinear error decays much faster than the theoretical
predicted rates and it can reach the strain increment $\Delta
\gamma=10^{-8}$.  This phenomenon has been observed in
\cite{OrtnerShapeev:2010} and is likely related to superconvergence of
local quantities of interest.  The indicator of the superconvergence
is the concentration of the critical eigenmode corresponding to
$\gamma^\c$ near the defect, which is illustrated in Figure
\ref{Fig:eigenmode}.

\begin{figure}
\centering
\includegraphics[width=10 cm ]{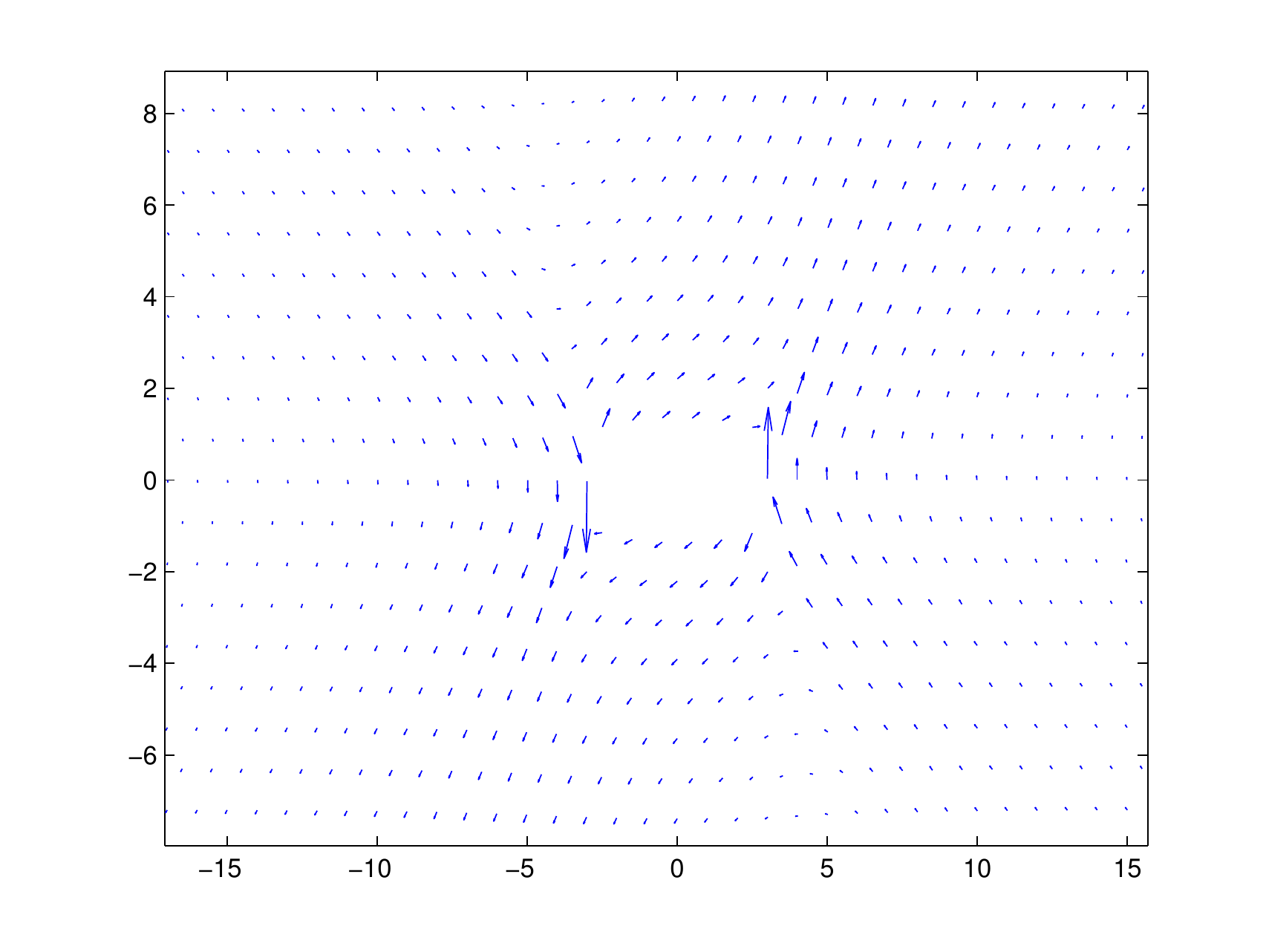}
\caption{The zoomed-in critical eigenvector of critical strain of vertically stretching a micro-crack.  We set $N=200$, crack length$=5$, strain increment
$\Delta\gamma=10^{-11}$ and $\alpha=4.$
}\label{Fig:eigenmode}
\end{figure}

We also study the relative errors of the critical strains for two
different choices of the blending width, $K=2$ and $K\approx
R_a^{3/5}+2$.  Motivated by the analysis in \cite{OrtnerShapeev:2010},
the size of the atomistic core is chosen to be $R_a=\sqrt{N}$.
According to Figure~\ref{Fig:CrackvsN}, the relative errors for
$K\approx R_a^{3/5}+2$ are approximately $10$ times smaller than those
for $K=2$. But both graphs decay rapidly as $N$ increases. The rate of
decay appears to be quadratic.
\begin{figure}[h]
\begin{center}
\includegraphics[width=11 cm ]{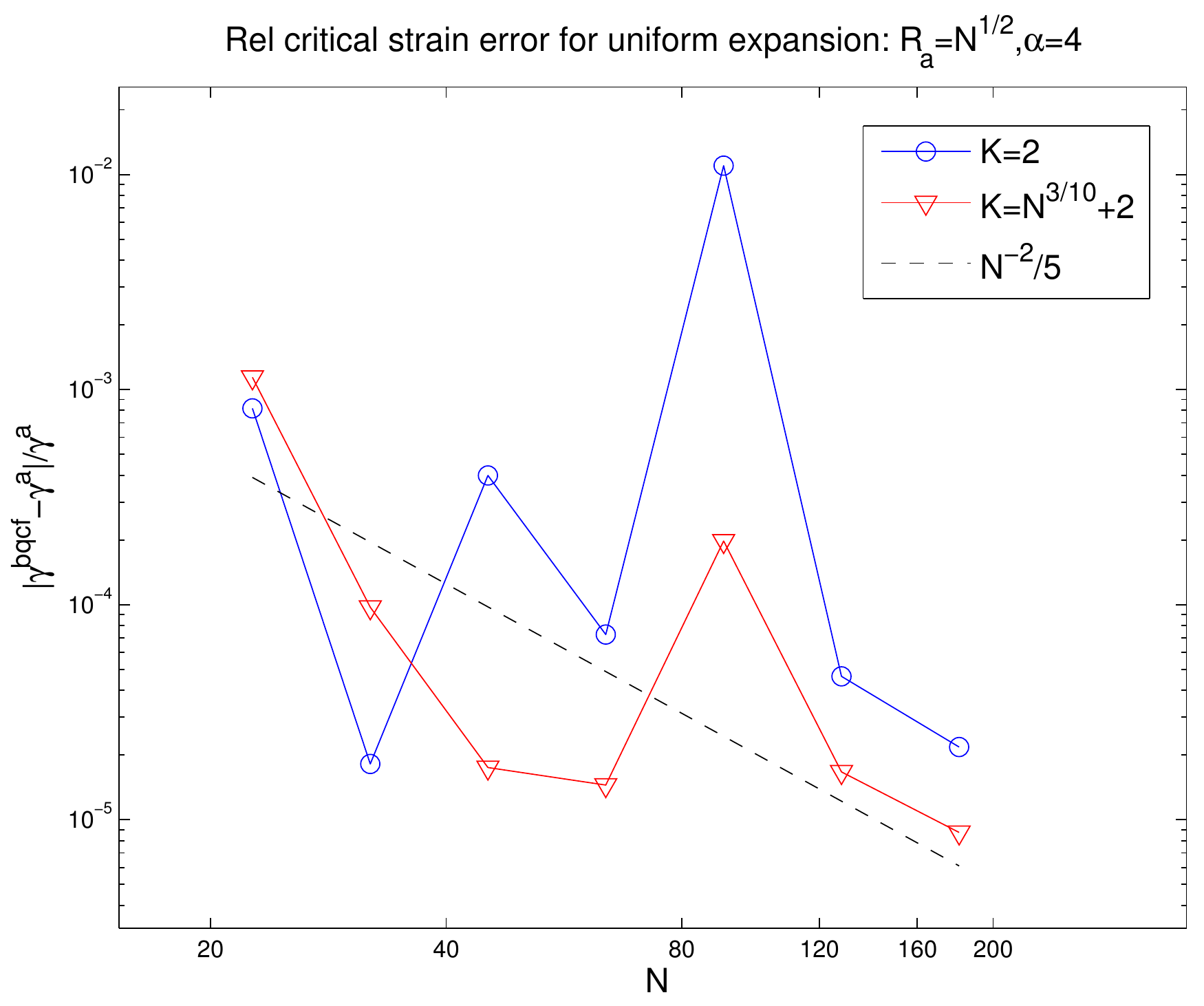}
\end{center}
\caption{The relative errors of the critical strains of vertically stretching a
  micro-crack in $\log_{10}$ scale plot.  We set crack length$=5$ {\helen{and
  $\alpha=4$}}.  $\gamma^{\rm a}$, $\gamma^{\rm{bqcf}}$ are the
  critical strains for the atomistic and B-QCF models,
  respectively. The dashed line is the theoretical
  asymptote.} \label{Fig:CrackvsN}
\end{figure}
\end{enumerate}

{\as{
\begin{remark}
The numerical computations in this section are conducted
without coarsening. The main reason for this was that
coarsening introduces more approximation parameters and
potentially more fluctuations in the results. Typically,
coarsening does not reduce the stability, and we therefore
expect that our stability results will remain valid for
\emph{any} coarsening. Another reason to discard coarsening was
to better compare the numerics with the theory. However, the
purpose of a/c coupling is to reduce the number of degrees of
freedom of an atomistic computation, therefore coarsening is
required when comparing efficiency (i.e., accuracy against the
number of degrees of freedom) of different methods.

\end{remark}
}
}
\def\Ea{\mathcal{E}^{\rm a}}
\def\Ec{\mathcal{E}^{\rm c}}
\def\L{\mathbb{L}}
\def\Vac{\mathcal{V}}
\def\Lv{\mathbb{L}_\Vac}
\def\Lom{\mathcal{L}}
\def\mA{B}
\def\R{\mathbb{R}}
\def\Th{\mathcal{T}_h}
\def\Nh{\mathcal{N}_h}
\def\Nhint{\Nh^{\rm free}}
\def\Ph{{\rm P}^1_h}
\def\Yh{{\rm Y}_h}

\section{The Accuracy of B-QCF}
\label{sec:accuracy}
In the previous sections, we investigated the {\em positivity} of the
B-QCF operator. One motivation for this study is that these
experiments fill a gap in our error analysis of the B-QCF method
\cite{VKOr:blend2}. We now briefly review these results and then
include some numerical experiments demonstrating the superior accuracy
of B-QCF over other a/c coupling schemes that we have investigated
previously in \cite{BQCEcomp}.

\subsection{Implementation of the B-QCF method}
\label{sec:bqcf_impl}
Let $\Vac\subset\L$ be a set of vacancy sites and $\Lv := \L \setminus
\Vac$ the corresponding lattice with defects. Let $\mA \in \R^{2
  \times 2}$ be the applied far-field strain. We consider the
atomistic problem
\begin{equation}
  \label{eq:atm_problem}
  y^{\rm a} \in \arg\min \big\{ \Ea(y) : y : \Lv \to \R^2, y(x) \sim B x \text{ as } |\xi|
  \to \infty \big\}.
\end{equation}
We remark that one must carefully renormalize $\Ea$ in order to
rigorously make sense of this problem; see e.g. \cite{VKOr:blend2} for
the details. \co{The vacancy sites are accounted for in the definition
  of $\Ea$ by simply removing the relevant pair interactions.}

We wish to approximate this problem with a practical variant \as{(i.e., with coarsening)} of the
B-QCF method. To that end, we choose $R_a, R_b = R_a + K, N \in
\mathbb{N}$ in such a way that all vacancy sites are contained in the
atomistic region $\Omega_a$, which is a hexagon with side length
$R_a$. The blending region is defined analogously. The full
computational domain is given by $\Omega$, which is a hexagon with
side length $N$. We triangulate $\Omega$ in such a way that it matches
the canonical triangulation of the triangular lattice in
$\Omega$.

Let $\Th$ denote the set of triangles, let $\Nh$ denote the nodes of
the triangulation, and let $\Nhint := \Nh \setminus (\Vac
\cup \partial\Omega)$ denote the {\em free nodes}.

Let $\Ph$ denote the space of all functions $v_h : \Omega \to \R^2$,
that are continuous and piecewise affine with respect to the
triangulation $\Th$. The space of admissible trial functions is then
given by
\begin{displaymath}
  \Yh := \big\{ y_h \in \Ph : y_h(x) = B x \text { for } x
  \in \partial\Omega \big\}.
\end{displaymath}
Each deformation $y_h \in \Yh$ is understood to be extended by $B x$
outside of $\Omega$ and thereby gives rise to an admissible atomistic
configuration.

We define the discretized Cauchy-Born energy functional as
\begin{displaymath}
  \Ec(y_h) := \sum_{T \in \Th} {\rm vol}(T) W\big(\nabla y_h|_T\big),
\end{displaymath}
where ${\rm vol}(T)$ in $2$D is the area of the triangle $T$.
We can define the discretized B-QCF operator, for a given blending
function $\beta$, as follows:
\begin{displaymath}
  F^{\rm bqcf}(x; y_h) :=
  (1-\beta(x)) \frac{\partial \Ea(y)}{\partial y(x)}\Big|_{y = y_h}
  + \beta(x) \frac{\partial \Ec(z_h)}{\partial z_h(x)}\Big|_{z_h =
    y_h} \qquad \text{ for } x \in \Nhint.
\end{displaymath}
In the B-QCF method, we aim to find a solution $y_h^{\rm bqcf} \in \Yh$ satisfying
\begin{equation}
  \label{eq:bqcf_method}
  F^{\rm bqcf}(x; y_h^{\rm bqcf}) = 0 \qquad \forall x \in \Nhint.
\end{equation}
We remark that this method has essentially five approximation
parameters that must be chosen carefully: the atomistic region size
$R_a$, the blending width $K$, the computational domain size $N$, the
blending function $\beta,$ and the finite element mesh $\Th$.

\subsubsection{Practical considerations}
\label{sec:implementation-more}
To implement \eqref{eq:bqcf_method} in practice, we need to specify
further details of the method:
\begin{enumerate}
\item In our choice of blending function, we deviate from the optimal
  choice of a $C^{2,1}$-blending function and instead choose only a
  $C^{1,1}$ blending function, which is more easily
  constructed. \co{This is justified, firstly, by our foregoing
    numerical experiments which suggest that little additional
    accuracy in the stability regions can be gained in the
    pre-asymptotic regime by using quintic splines (i.e.,
    $C^{2,1}$-blending), and secondly, because the consistency error
    does not depend on the regularity of the blending function.}

  We choose the blending function proposed in \cite{BQCEcomp}, which
  minimizes $\| \nabla^2 \beta \|_{L^2}$, or a discrete variant
  thereof, in a precomputation step (see \cite{BQCEcomp} for the
  details).


\item In addition to the blending region $\Omega_b$ we ensure that two
  additional ``layers'' of atoms outside of it belong to $\Nh$. This
  makes the implementation of the atomistic force contribution in
  \eqref{eq:bqcf_method} straightforward.

  Moreover, we ensure that the vacancy sites do not affect the forces
  on atoms $x$ where $\beta(x) \neq 0$. This ensures that all the
  Cauchy-Born force contributions in \eqref{eq:bqcf_method} are
  the correct Cauchy-Born forces.

\item To obtain an appropriate initial guess for the B-QCF solutions,
  we first solve the corresponding energy-based blended QCE method
  (B-QCE)~\cite{BQCEcomp} with the same approximation parameters,
  using a preconditioned line search method. The details are described
  in \cite{BQCEcomp}. The B-QCE solution is then taken as a starting
  guess for the B-QCF Newton iteration to solve
  \eqref{eq:bqcf_method}.
\co{If no B-QCE code is readily available, then a natural
    alternative would be to implement a damped Newton method for
    B-QCF. }

  We remark, that the Jacobian matrix of the B-QCF operator is
  straightforward to assemble from the Hessians of the atomistic and
  Cauchy-Born energy. Nevertheless, for large 3D simulations, more
  sophisticated solution methods may be required.

\item We are now only left to choose the remaining approximation
  parameters $R_a, K, N$ and the mesh $\Th$.
\end{enumerate}

\subsection{Error versus computational cost}
\label{sec:errest}
We briefly review the main ideas of our analysis in \cite{VKOr:blend2}
without technical details. A first key result is that if the atomistic
solution is stable ($\delta^2\Ea(y^{\rm a})$ is positive
definite) and the linearized B-QCF operator $\delta F^{\rm bqcf}(\cdot; B
x)$ is positive definite, then choosing $R_a$ and $K$ sufficiently large
implies that $\delta F^{\rm bqcf}(\cdot; y^{\rm a})$ is also positive
definite, that is, the B-QCF method is {\em stable} under these
conditions. To achieve this in practice, we need to choose $K^3 \gg
R_a$
(recall from section \ref{sec:implementation-more} that we have chosen
a sub-optimal $\beta$).

From this stability result, we can deduce the existence of a B-QCF
solution in a neighborhood of the atomistic solution, and an error
estimate in terms of the {\em best approximation error} (the best
approximation of $y^{\rm a}$ from the finite element space $\Yh$). and of the
modeling error (the force discrepancy of the B-QCF and atomistic
models). We estimate the error in the strain $\nabla y^{\rm a} -
\nabla y_h^{\rm bqcf}$ in terms of the ``smoothness'' of $y^{\rm a}$,
which is measured in terms of bounds on the derivatives $\nabla^j
y^{\rm a}$. The derivatives of the discrete functions $y^{\rm a}$ are
understood as derivatives of a smooth interpolant. (See
\cite{VKOr:blend2} for the details.)

Dropping an unimportant term for the sake of readability, our error
estimate reads
\begin{align}\label{error.bqcf}
  \| \nabla y^{\rm a} - \nabla y_h^{\rm bqcf} \|_{L^2(\R^2)} \lesssim
  C^{\rm stab} \Big( C^\beta \| \nabla^3 y^{\rm a}
  \|_{L^2(\R^2 \setminus \omega_{\rm a})} + \| h \nabla^2 y^{\rm a}
  \|_{L^2(\Omega \setminus \omega_{\rm
      a})} + \| \nabla y^{\rm a} \|_{L^2(\R^2 \setminus \omega)} \Big),
\end{align}
where $\| \nabla^3 y^{\rm a} \|_{L^2(\R^2 \setminus \omega_{\rm a})}$
measures the modeling error, $\| h \nabla^2 y^{\rm a} \|_{L^2(\Omega
  \setminus \omega_{\rm a})}$ the finite element discretization error
and $\| \nabla y^{\rm a} \|_{L^2(\R^2 \setminus \omega)}$ the error in the
far-field due to the artificial boundary condition (the two latter
errors comprise the best approximation error). The domains
$\omega_{\rm a}, \omega$ are slightly smaller hexagonal subsets of,
respectively, $\Omega_{\rm a}$ and $\Omega$, with comparable
side lengths.

In addition, $C^{\rm stab}$ is a stability constant that is uniformly
bounded for $R_a
\ll K^3$, and
\begin{displaymath}
  C^\beta := K^{-1/2} R_a^{1/2} \log\big|R_a / N\big|
\end{displaymath}
is a $\beta$-dependent prefactor, which arises from a crucial
inequality, $\| \nabla(\beta v) \|_{L^2} \leq C^\beta \| \nabla v
\|_{L^2}$, in the consistency analysis of B-QCF.

We choose $K \approx R_a$ and $N$ a
polynomial of $R_a$ (we will see momentarily why this is natural),
then $C^\beta$ is uniformly bounded and in addition, we choose $R_a
\ll K^3$
, which we require for stability. With this choice, it is easy
to see that $C^\beta \| \nabla^3 y^{\rm a} \|_{L^2(\R^2 \setminus
  \omega_{\rm a})} \lesssim \| h \nabla^2 y^{\rm a} \|_{L^2(\Omega
  \setminus \omega_{\rm a})}$ (recall that we are working in units
where atomic spacing is $1$), and hence we can simply ignore the
modeling error term from now on.

We recall from \cite{BQCEcomp} that the atomistic method (ATM) is given by the B-QCF method
with $\beta\equiv 0.$ We also recall the
corresponding error estimates (dropping less important terms) for the atomistic
(ATM) and the B-QCE methods~\cite{VKOr:blend2,BQCEcomp}
\begin{align}\label{error.bqce}
 \| \nabla y^{\rm a} - \nabla y^{\rm atm} \|_{L^2(\R^2)}&\lesssim
  \| \nabla y^{\rm a} \|_{L^2(\R^2 \setminus \omega)} ,\\ \label{error.atm}
 \| \nabla y^{\rm a} - \nabla y_h^{\rm bqce} \|_{L^2(\R^2)}&\lesssim
   \| \nabla^2 \beta
  \|_{L^2(\R^2 \setminus \omega_{\rm a})} + \| h \nabla^2 y^{\rm a}
  \|_{L^2(\Omega \setminus \omega_{\rm
      a})} + \| \nabla y^{\rm a} \|_{L^2(\R^2 \setminus \omega)} .
\end{align}

To better understand the best approximation error, we need to understand the
regularity of $y^{\rm a}$. Since the problems only involve defects
with zero Burgers vector, it is reasonable to assume based on linear
elasticity, that
\begin{displaymath}
  |\nabla^j y^{\rm a}(x)| \sim |x|^{-j-1}.
\end{displaymath}
\co{(We stress that this estimate only applies in the far-field. In
  the preasymptotic regime different rates of decay might be observed,
  e.g., $|\nabla^j y^{\rm a}(x)| \sim |x|^{1/2-j}$ for the micro-crack
  case discussed in~\S~\ref{sec:microcrack}.)}

Having this explicit knowledge about the elastic field, we can
optimize our choice of finite element triangulation. Using the
construction in \cite{OrtnerShapeev:2010} and also used successfully in
our B-QCE experiments in \cite{BQCEcomp}, we obtain a triangulation
$\Th$ (as a function of $R_b$ and $N$), for which the following
estimate holds:
\begin{displaymath}
  \| h \nabla^2 y^{\rm a} \|_{L^2(\Omega
  \setminus \omega_{\rm a})} + \| \nabla y^{\rm a} \|_{L^2(\R^2
  \setminus \omega)} \lesssim R_a^{-2} + N^{-1}.
\end{displaymath}
Thus, we choose $N \approx R_a^2$ to balance these two error
contributions.

Finally, we note that, with this construction, the number of degrees
of freedom in $\Yh$, ${\rm DoF} := {\rm dim} \Yh = 2\#\Nhint$ is
approximately equal to ${\rm DoF} \approx R_a^2$. (In particular, the
number of degrees of freedom in the atomistic, blending and continuum
regions are comparable.)

In summary, choosing $K \approx R_a, N \approx R_a^2$, the blending
function $\beta$ according to the construction proposed in
\cite{BQCEcomp}, and the finite element mesh according to the
construction proposed in \cite{OrtnerShapeev:2010}, we obtain from \eqref{error.bqcf} the
error estimate
\begin{equation}
  \label{eq:errest}
  \| \nabla y^{\rm a} - \nabla y_h^{\rm bqcf} \|_{L^2(\R^2)} \lesssim
   {\rm DoF}^{-1}.
\end{equation}
We note that $N=R_a$ in the ATM method, and consequently we obtain from
\eqref{error.atm}
\begin{equation}
  \label{eq:errest2}
  \| \nabla y^{\rm a} - \nabla y^{\rm atm} \|_{L^2(\R^2)} \lesssim
   {\rm DoF}^{-1/2};
\end{equation}
thus demonstrating an improved rate of convergence for the B-QCF method
in comparison with the ATM method.

We remark that this is optimal for P1-finite element type
coarse-graining schemes, as the modeling error is in fact dominated
by the finite element error. In particular, it is a substantial
improvement over the B-QCE method, for which the corresponding error
estimate obtained from \eqref{error.bqce} is
\[
\| \nabla y^{\rm a} - \nabla y_h^{\rm bqce}
\|_{L^2(\R^2)} \lesssim {\rm DoF}^{-1/2}.
\]
We note that the B-QCE method can be shown to have a higher rate of convergence
than the ATM method for defects with nonzero Burgers vector (such as dislocations) which
have a lower rate of decay.  The finite element coarse-graining of the B-QCE method
can more efficiently approximate the larger region where the strain gradient is
significant; \cite{VKOr:blend2,BQCEcomp} for the details.

\subsection{Numerical rates}
We test our analytical predictions against the two numerical examples,
for which we already tested the B-QCE method in \cite{BQCEcomp}. In
both examples, we choose the Morse interaction potential
{\helen{
\begin{displaymath}
  \phi(r) =\left[1-\exp(-\alpha(r-1))\right]^2,
\end{displaymath}
}}
with stiffness parameter $\alpha = 4$.

We compare the B-QCF method with a pure atomistic computation on a
finite domain, with the QCE and B-QCE methods (cf. \cite{BQCEcomp} for
a detailed description of these three methods) and with the pure QCF
method, which is simply the B-QCF method with $K < 1$ (i.e., $\beta(x)
\in \{0, 1\}$).

Finally, we have also included a highly optimized B-QCE variant where
we choose $K \approx R_a^2$ and $N \approx R_a^4$, which is a very
unexpected scaling, but yields improved errors in the preasymptotic
regime; see \cite[Remark 4.3]{BQCEcomp}. We denote this method by
B-QCE+ in the error graphs.

\subsubsection{The di-vacancy example}
\label{sec:divacancy}
\begin{figure}
  \begin{center}
    \includegraphics[height=6cm]{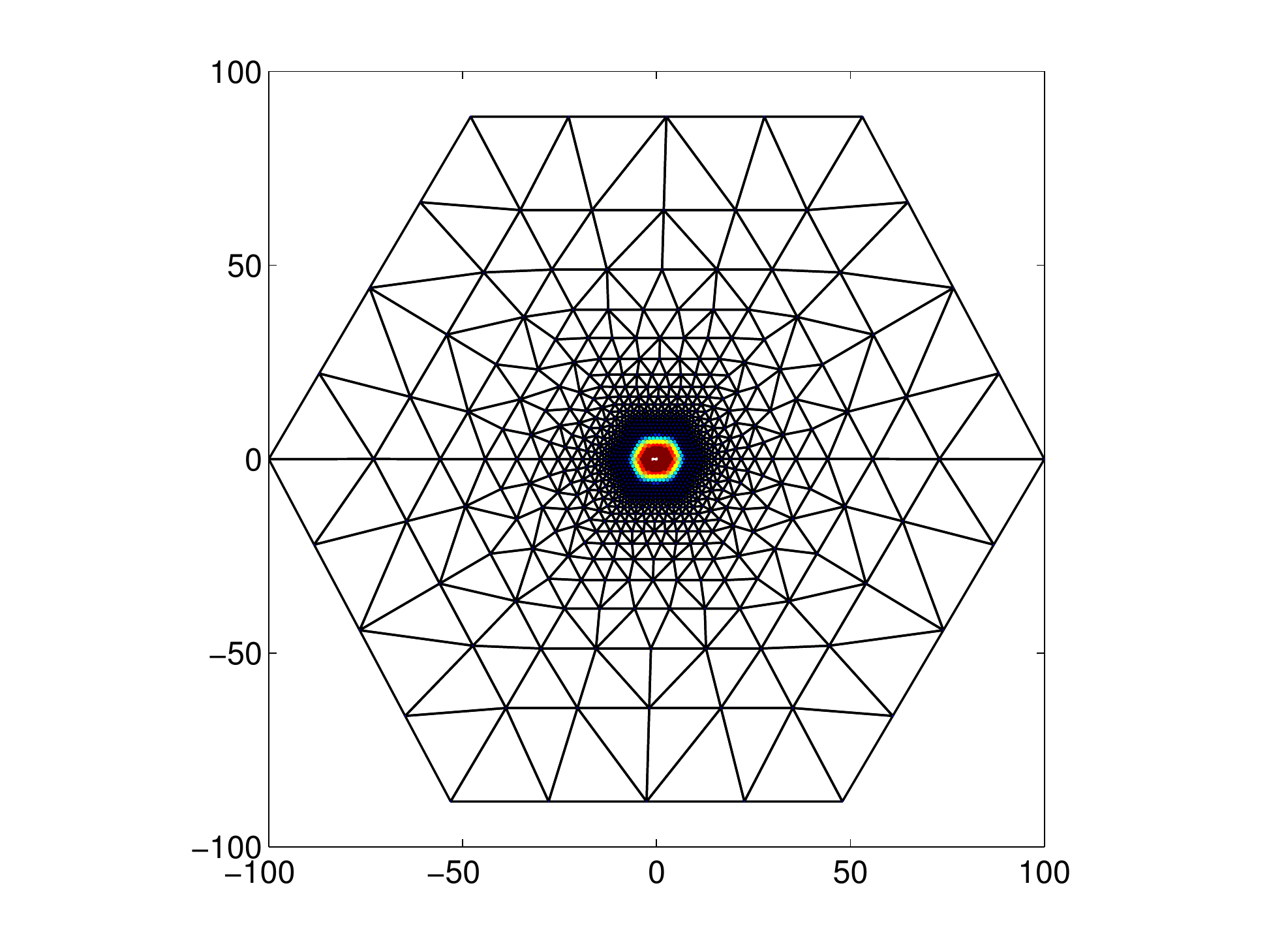} \quad
    \includegraphics[height=6cm]{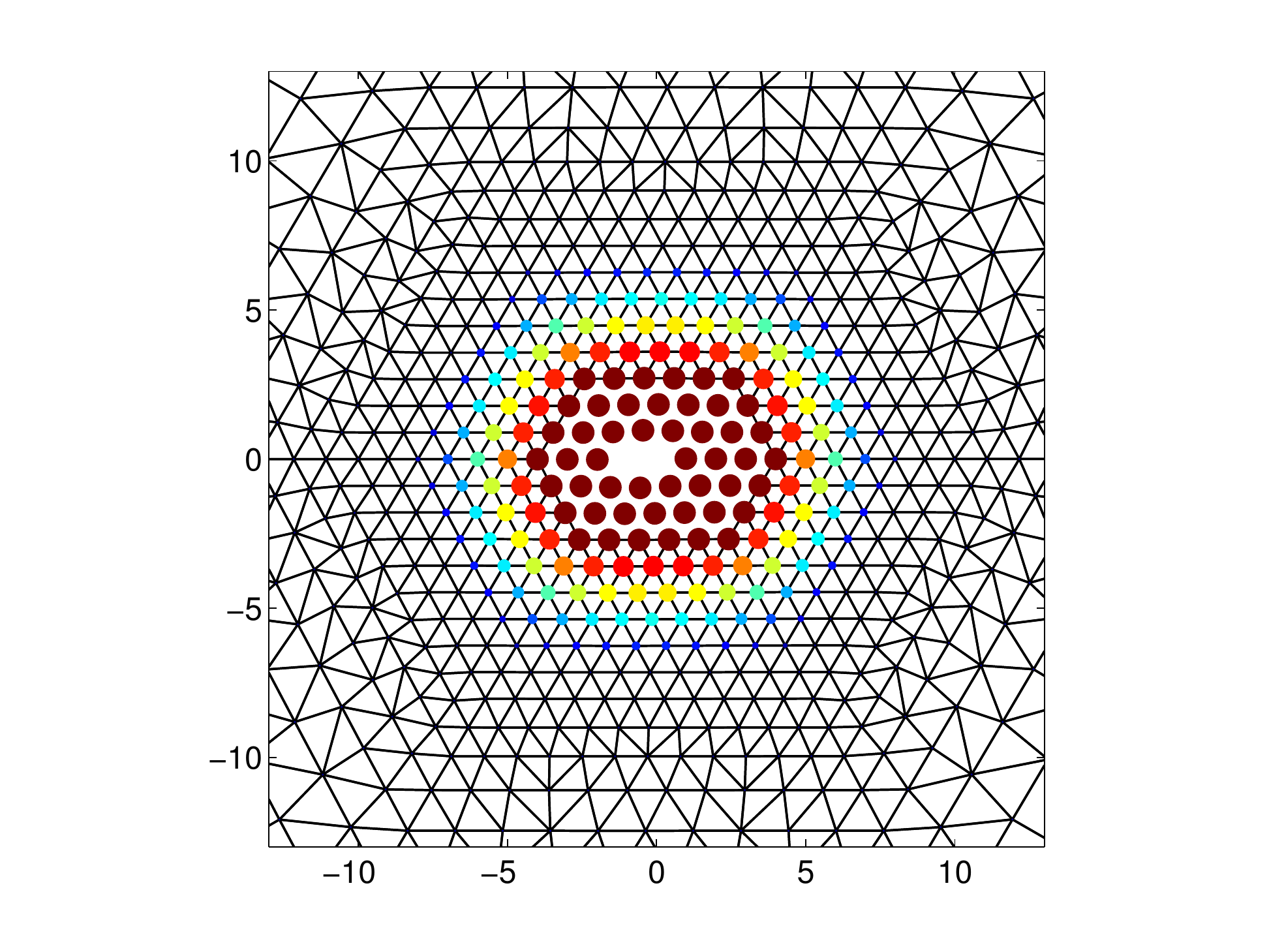}
    \caption{\label{fig:divac_setup} Setup of the B-QCF method for the
      di-vacancy example, for a specific choice of approximation
      parameters, shown in deformed equilibrium. The size/color of the
      atoms in the center correspond to decreasing values of
      $(1-\beta(x))$. }
  \end{center}
\end{figure}
We choose the vacancy set $\Vac = \{0, e_1\}$ and the macroscopic
strain
\begin{displaymath}
  \mA = \left(\begin{matrix}
    1.03 &  0.3 \\
    0.0 & 1.03
  \end{matrix}\right) \cdot
  \mA_0,
\end{displaymath}
where $\mA_0$ is a minimizer of $W$ ($3\%$ uniform
stretch and $3\%$ shear from ground state). The setup of the B-QCF method for the
di-vacancy problem is shown in Figure \ref{fig:divac_setup}.

\begin{figure}
  \begin{center}
    \includegraphics[width=10cm]{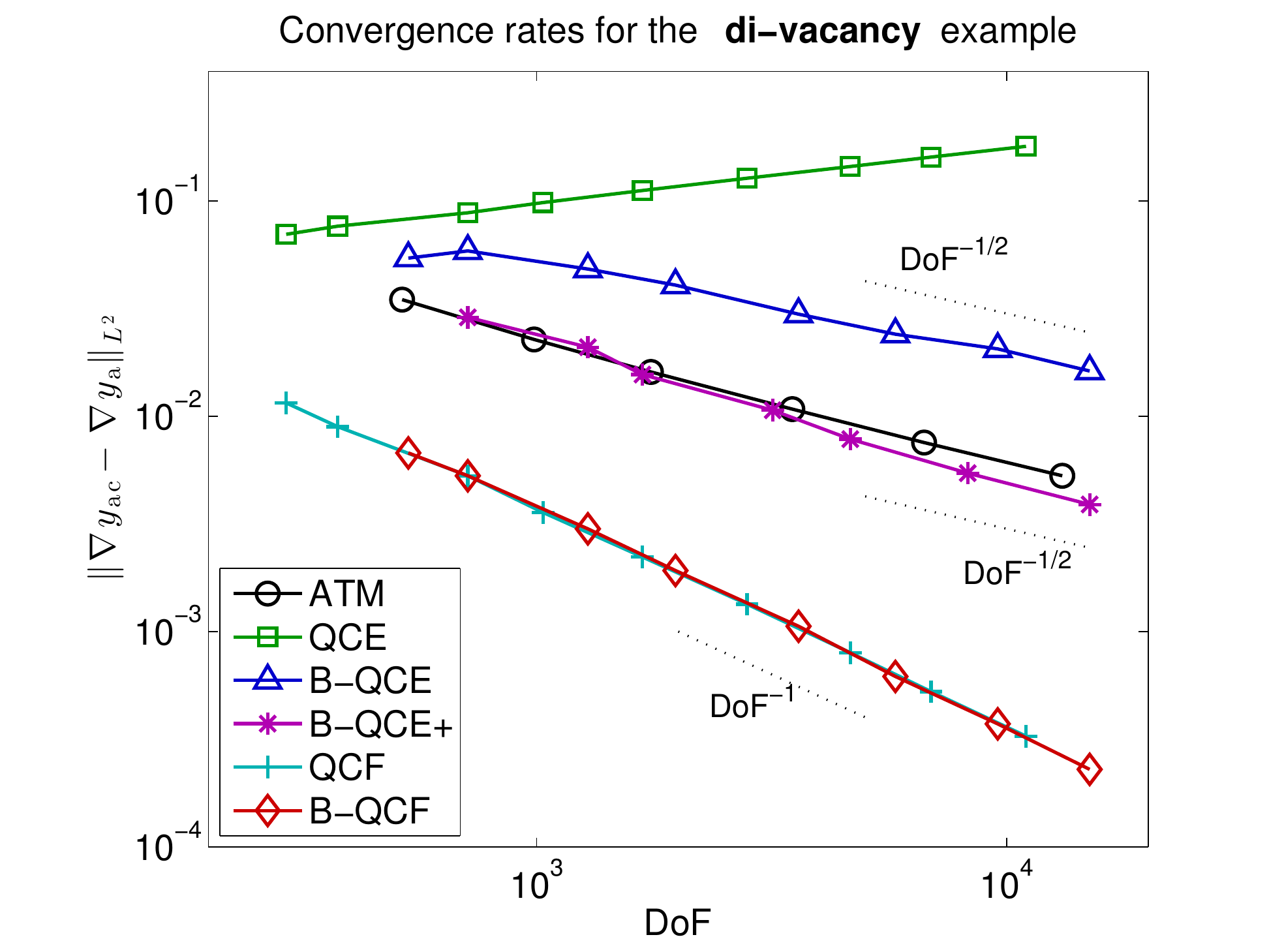}
    \caption{\label{fig:divac_err2} Plots of computational cost (DoF)
      versus error in the energy-norm for various a/c coupling methods
      approximating the di-vacancy problem described in
      section \ref{sec:divacancy}.}
  \end{center}
\end{figure}

In Figure \ref{fig:divac_err2}, we plot the degrees of freedom (DoF)
against the error in the energy-norm for the various a/c coupling
methods that we consider. As predicted by our analysis, the B-QCF
method clearly outperforms all other methods, with the exception of
the QCF method, which is barely distinguishable from the B-QCF method
in this graph. Unfortunately, we cannot offer a satisfactory theory
for the QCF method at present.

We also remark that, due to the high consistency error committed in
the interface region, the B-QCE does not even outperform a plain
atomistic computation in this particular example. (But it will clearly
outperform the fully atomistic method (ATM) in the micro-crack example, where the
elastic field is much more significant.)

\subsubsection{The micro-crack example}
\label{sec:microcrack}
\begin{figure}
  \begin{center}
    \includegraphics[height=6cm]{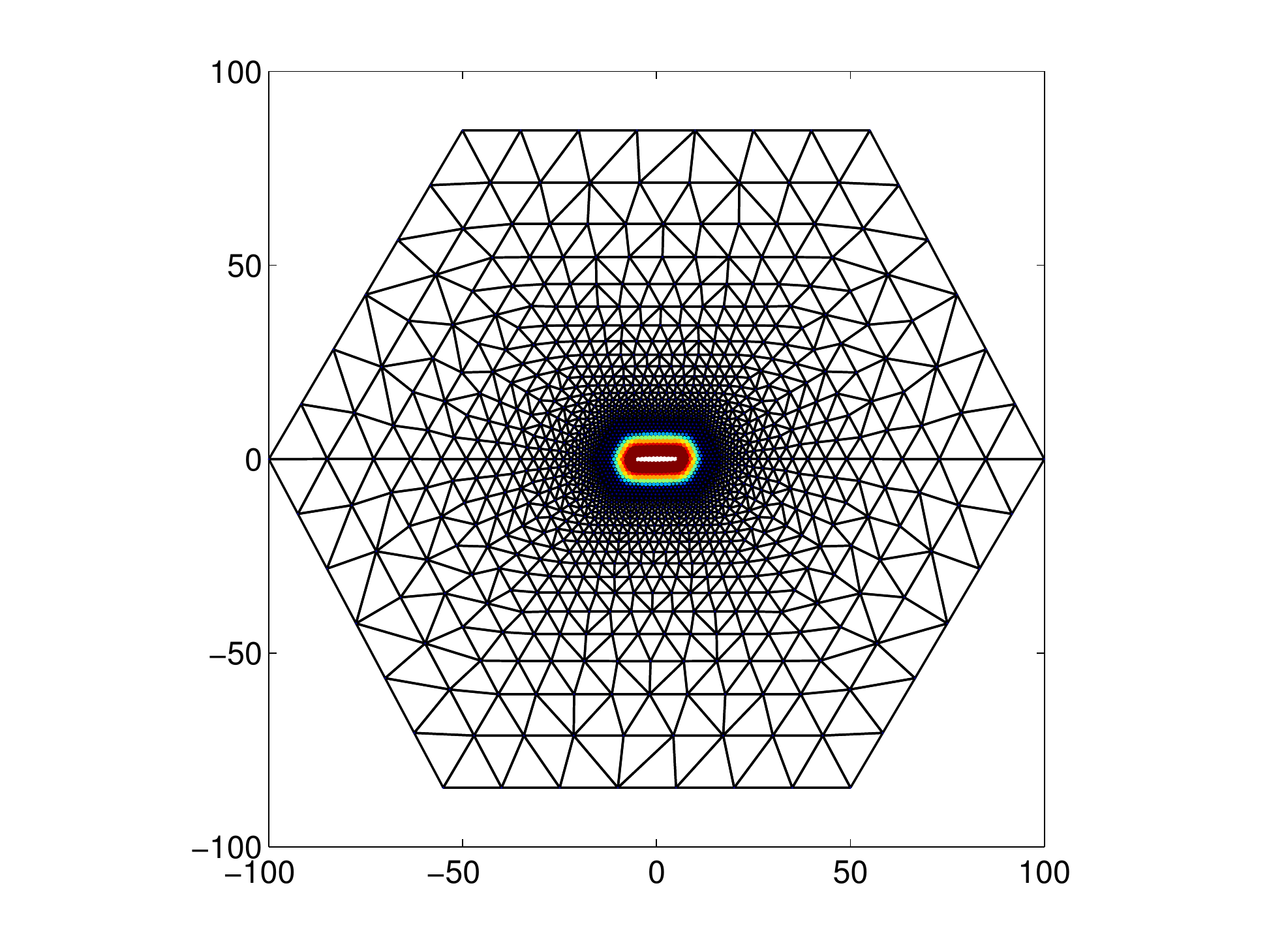} \quad
    \includegraphics[height=6cm]{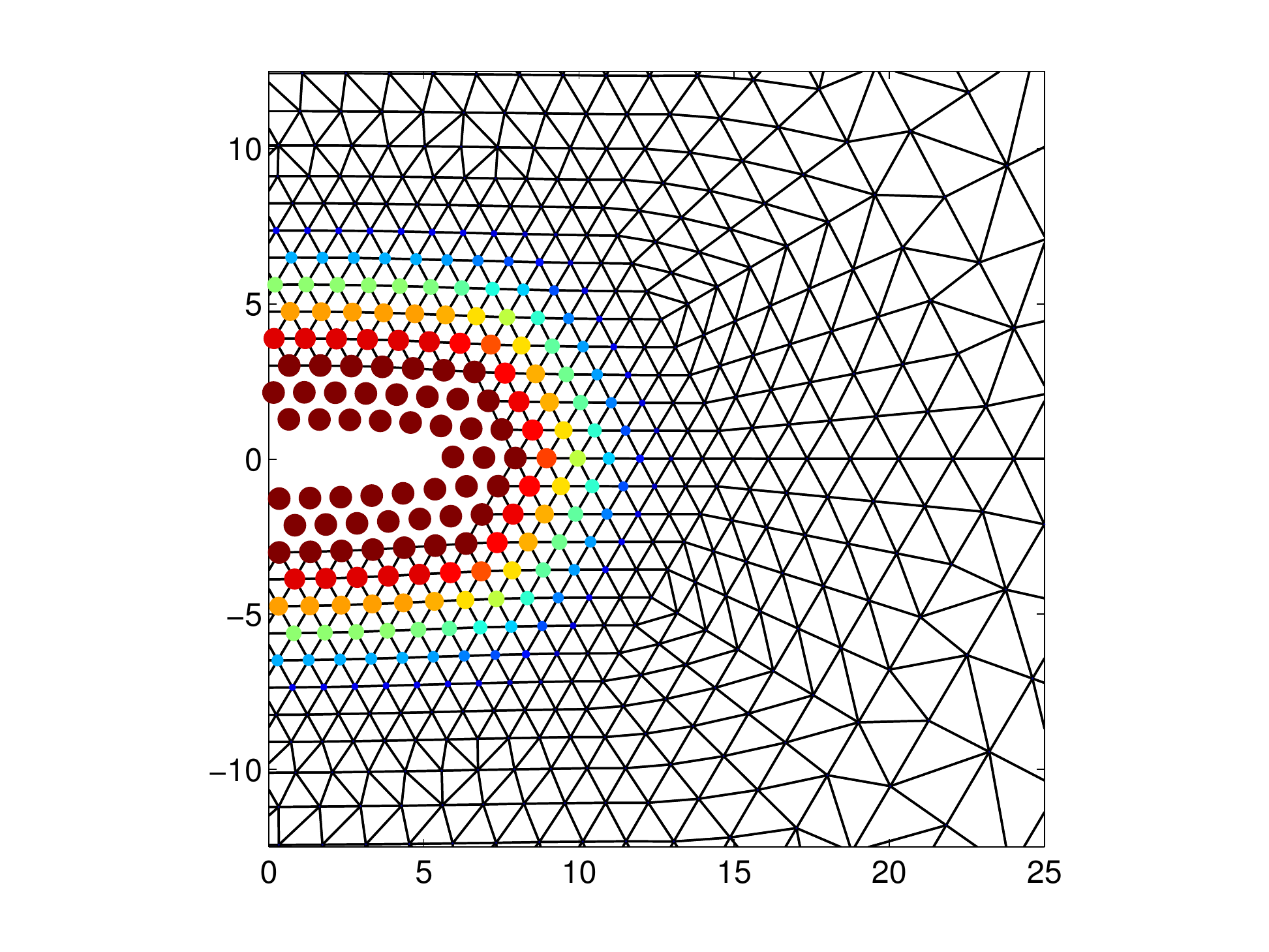}
    \caption{\label{fig:mcrack_setup} Setup of the B-QCF method for
      the micro-crack example, for a specific choice of approximation
      parameters, shown in deformed equilibrium. The size/color of the
      atoms in the center correspond to decreasing values of
      $(1-\beta(x))$.}
  \end{center}
\end{figure}
In the micro-crack \co{(or void)} example, we choose the
vacancy set $\Vac = \{-5 e_1, \dots, 5 e_1\}$ and the
macroscopic strain
\begin{displaymath}
  \mA = \left(\begin{matrix}
    1.0 &  0.03 \\
    0.0 & 1.03
  \end{matrix}\right) \cdot
  \mA_0,
\end{displaymath}
where $\mA_0$ is a minimizer of $W$ ($3\%$ tensile stretch and $3\%$
shear from ground state). The setup of the B-QCF method for the
micro-crack problem is shown in Figure \ref{fig:mcrack_setup}.

\begin{figure}
  \begin{center}
    \includegraphics[width=10cm]{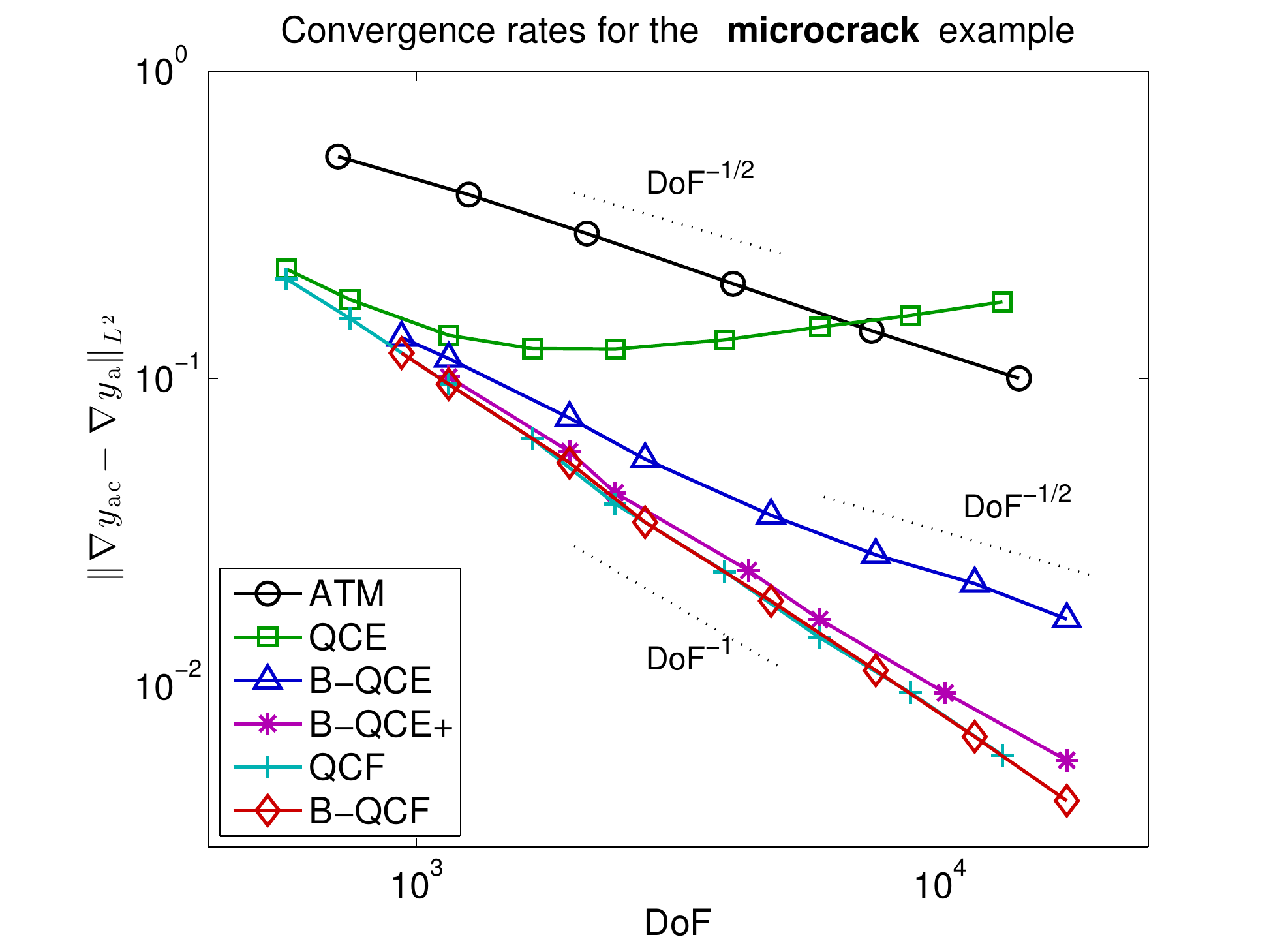}
    \caption{\label{fig:mcrack_err2} Plots of computational cost (DoF)
      versus error in the energy-norm for various a/c coupling methods
      approximating the micro-crack problem described in
      section \ref{sec:microcrack}.}
  \end{center}
\end{figure}

In Figure \ref{fig:mcrack_err2} we plot the degrees of freedom (DoF)
against the error in the energy-norm, for the various a/c coupling
methods that we consider. In this example the picture is less clear
than in the di-vacancy example due to a more significant preasymptotic
regime, which is caused by the more significant deformation admitted
by the microcrack. In the preasymptotic regime we observe that the QCE
and B-QCE methods perform much better than expected, but eventually
fall back to the predicted rates. By contrast, the B-QCF and QCF
methods display clear systematic convergence at the predicted rate
throughout.

We also note that, in this example, the B-QCE+ method performs
comparable to the B-QCF and QCF methods, at least in the preasymptotic
regime accessible in the experiment.

\section{Conclusion}
We have formulated an atomistic-to-continuum force-based coupling, which we call the
blended force-based quasicontinuum (B-QCF) method. In this paper,
we numerically studied the stability as well as accuracy of the B-QCF method.
We computed the
critical strain errors between the atomistic and B-QCF models with different sizes of
the blending region under different types of deformations.

The main theoretical conclusion in \cite{BQCF} is that the required blending width
to ensure coercivity of the linearized B-QCF operator is surprisingly
small. For both $1$D and $2$D uniform expansion, the
computational results of the linearized operators perfectly
match the analytic predictions. In addition, the stability for a general class of
homogeneous deformations of the $2$D B-QCF operator becomes almost the same as that of
the atomistic model by using a very small blending region, in contrast to the fact that
the stability region of the force-based quasicontinuum (QCF) method, that is, the B-QCF method without blending region,
is just a proper subset of the fully atomistic model.
However, the critical strain error for the B-QCF operator applied to shear deformation
seems to only linearly depend on the system size and is thus insensitive to blending width.

For the problem of a microcrack in
a two-dimensional crystal, we studied the nonlinear stability of the B-QCF operators.
The critical strain error decays faster than the prediction,
and it can be as small as the strain increment. However, we find that the error increases a little bit
when the blending size becomes larger, which is possibly due to round-off error.

Moreover, we implemented a {\em practical} version of the B-QCF method.
We briefly reviewed the accuracy results in terms of computational cost \cite{VKOr:blend2}.
The numerical experiments, di-vacancy and microcrack demonstrate the superior accuracy
of B-QCF over other a/c coupling schemes that we have investigated
previously in \cite{BQCEcomp}.

The BQCF method with a surprisingly
small blending region is an appealing choice for numerical simulations of atomistic multi-scale problems
as it is always consistent and can be guaranteed by both theory and benchmark testing
 to be positive definite when the fully atomistic operator
is positive definite.


\section{Acknowledgments}
We appreciate helpful discussions with Brian Van Koten.

\bibliographystyle{abbrv}
\bibliography{BQCF_cmame_rev6}
\end{document}